
\input amstex
\documentstyle{amsppt}

\define\cp{\overline{\bold{CP}}^2}
\define\Fr{Fr\wedge^2_+TX}
\define\Fl{\tilde\wedge^2_+}
\define\Mod #1 #2{{\Cal M}(P_{#1},g_{#2})}
\define\FMod #1 #2 #3{{\Cal M}^{#3}(P_{#1},g_{#2})}
\define\CMod #1#2{{\overline{\Cal M}}(P_{#1},g_{#2})}
\define\HMod #1#2{W(P_{#1},g_{#2})}
\define\FHMod #1#2#3{W^{#3}(P_{#1},g_{#2})}
\define\mubar{\overline\mu}
\define\bZ {\Bbb Z}
\define\bQ {\Bbb Q}
\define\bR {\Bbb R}
\define\bC {\Bbb C}
\define\bH {\Bbb H}
\define\bP {\Bbb P}
\define\bE {\Bbb E}
\define \br {\overline}

\topmatter
\title Donaldson Wall-Crossing Formulas Via Topology
\endtitle
\author Thomas Leness \endauthor
\affil Michigan State University \endaffil
\address East Lansing, MI 48824-1027\endaddress
\email leness\@math.msu.edu \endemail
\date March 26th, 1996 \enddate
\abstract The wall-crossing formula for Donaldson invariants of
simply conected four manifolds with $b^+=1$ is shown to be
a topological invariant of the manifold, for reducibles with
two or fewer singular points. The explicit formulas derived agree
with those of \cite{FQ} and \cite{EG}.\endabstract
\endtopmatter
\document
\hoffset=0.5cm
\hsize=6.5in
\noindent
{\bf 0: Definitions and Notation}

\

Although the Seiberg-Witten invariants have replaced Donaldson invariants
as the tool of choice for answering questions about the differential
topology of four manifolds, the Donaldson invariants of simply
connected, smooth four manifolds with $b^+=1$, as
defined in \cite{Ko},\cite{KoM}, still pose interesting
questions.  These manifolds are the only known examples of
manifolds whose Donaldson invariants are not of simple type (see \cite{KoL}).
Unlike manifolds with $b^+>1$, the Donaldson invariants
of manifolds with $b^+=1$ are not metric independent, but vary with the
metric according to a \lq\lq wall-crossing formula\rq\rq.  This wall-crossing
formula is of interest in its own right, especially in
comparing the Donaldson and Seiberg-Witten moduli spaces.

Let $P\rightarrow X$ be an $SO(3)$ principal bundle with $p_1(P)=p$
over a smooth, simply connected four manifold $X$, $b^2_+(X)=1$,
and $c$ an integer lift
of $w_2(P)$.  $\CMod  p X$ will denote the Uhlenbeck compactification of
the moduli space of connections on $P$ which are anti-self-dual with respect to
the
metric $g_X$.  If there are no reducible connections in $\Mod p X$, then
the map $\mu : H_2(X;\bZ )\rightarrow H^2({\Mod p X};\bQ )$ is defined
by taking the slant product of the Pontrjagin class of the universal
bundle ${\Bbb P}\rightarrow {\Mod p X}\times X$, $\mu(x)=-{1\over 4}p_1({\Bbb
P})/x$.
If there are reducibles, the space giving the universal bundle is no longer a
bundle.  It is shown in \cite {DK, Thm 5.1.21} that the $\mu(x)$ class can be
non-trivial on the link of a reducible and thus not be extendable across the
reducible.
For a generic metric $g_X$, there
will be no reducibles and the cohomology class $\mu(x)$ extends to
$\mubar (x)\in H^2({\CMod p X};\bQ)$ as described in \cite {FM}.
One would like to define Donaldson invariants by pairing these classes
$\mubar (x)$ with the homology class given by the compactified moduli space.
This pairing will not be independent of the metric: if one tries to use the
moduli space associated to a path of metrics $g_t$ to form a cobordism between
the moduli spaces of the endpoint metrics,
$$\CMod p \lambda =\{([A],t): [A]\in \CMod p t \}$$
there may be reducibles in this space, so that the $\mubar(x)$ class would not
be
defined on it.

A reducible connection gives a reduction of $P$ to an
$S^1$ bundle $Q_\alpha$, $c_1(Q_\alpha)=\alpha$, $P=Q_\alpha\times_{S^1}SO(3)$.
The curvature of the reducible connection
must then be an ASD harmonic two form, thus perpendicular to the
ray of self-dual two forms $\omega(g_X)$ which are harmonic with respect to the
metric
$g$.  There will be a reducible connection $[\alpha]$ corresponding to
$\pm\alpha$ in $\CMod p X$
only if
$$\omega(g_X)\in W^\alpha=\{x\in H^2(X;\bR): x^2>0, x\cdot \alpha=0\}$$
where $H^2(X;\bR)$ is identified with the space of harmonic forms.
We call $W^\alpha$ the wall associated to $\alpha$.

A reducible connection $[\alpha]$ need not be in the highest stratum of
$\CMod P X$.  If $p=p_1(P)$, and we write $r=(\alpha^2-p)/4$, then
the reducible $[\alpha]$ will lie in the strata of
$\CMod P X$ contained in ${\Cal M}(P',g)\times \Sigma^r(X)$,
where $p_1(P')=p+4r$ and $w_2(P')=w_2(P)$.
The set of walls $W^\alpha$ where $\alpha\thinspace\text{mod}\thinspace
2=w_2(P)$
and $\alpha^2> p$ is called the set of $P$-walls.
The set of chambers $\Cal C$ of $X$ is the set of components of the complement,
in the positive cone of $H^2(X;\bR )$, of the $P$-walls.  The Donaldson
invariant associated
to the bundle $P$ will change only when the period point of the metric passes
through a $P$-wall.  We write $\tilde \Cal C$ for the pre-image of $\Cal C$
in the space of metrics
under the period map, so the Donaldson invariant depends only on the
element of $\tilde \Cal C$ in which the metric lives.
We can then define the Donaldson invariants
$$D^X_{d,c}: {\tilde \Cal C}\rightarrow \text{Sym}^d(H_2(X;\bZ)\oplus
H_0(X;\bZ))$$
to be
$$D^X_{d,c}(\tilde C)(z)=<\mubar (z),[\CMod p X]>.$$
Here $g_X$ is a metric such that $g\in\tilde C$, $d=-p+3$.  The integer lift of
$w_2(P)$, $c$ and a choice of an orientation of $H_0(X)\otimes H^2_+(X)$ give
an
orientation to the moduli space.

A reducible anti-self dual connection $[\alpha]$
corresponding to the line bundle $Q_\alpha$ with $c_1(Q_\alpha )=\alpha$
where $\alpha^2-p=4r$ gives a family of reducibles:
$$[\alpha ]\times \Sigma^r(X)$$
in the compactification.  Let $D(\alpha,g_0)$ denote the link of this
family of reducibles in the parametrized moduli space $\CMod p \lambda$,
then we define the difference term $\delta_P(\alpha,g_0)$ by
$$\delta_P(\alpha,g_0)(z)= < \mubar(z),[D(\alpha,g_0)]>.$$
The change in the Donaldson invariant along a path of metrics $g_t$
is given by
$$D^X_{d,c}(g_1)-D^X_{d,c}(g_{-1})=\sum_{\alpha}\epsilon(c,\alpha)\delta(\alpha,g_\alpha)$$
where $g_\alpha$ is the metric in the path $g_t$ where the reducible $\alpha$
occurs and $\epsilon(c,\alpha)=(-1)^{((c-\alpha)^2/2)}$.  This sum is over all
reducibles that appear in the path of metrics.

In \cite{KoM}, an argument is given that the difference term
$\delta_P(\alpha,g)$
is independent of the metric $g$ and depends only on the reducible.
This is necessary to to show that the Donaldson invariant depends not on
the metric component $\tilde C$ but only the chamber $C$.
There are
some technical difficulties with the argument given there
(these difficulties are outlined in Section 3.4).  Although the result
is almost certainly true, we do not use it in this calculation.

The correction terms $\delta_P(\alpha,g)$ remain mysterious, because of the
difficulty
in describing neighborhoods of a reducible.  The lower the stratum in which the
reducible appears, the harder the calculation appears to be.  That is, if
$\alpha^2-p_1=4r$, then $r$, the number of points of totally concentrated
curvature,
measures the complexity.
For $r=0$, the difference term has been calculated
in \cite{Do,Ko},for $r=1$ by Yang in \cite{Y} for $SU(2)$ bundles.
Calculations for
larger $r$ in the case of algebraic manifolds have
been done in \cite{EG} and \cite{FQ}.
In this paper, an expression for reducibles with $r=1,2$ is obtained, matching
those
in \cite{FR} and \cite{EG}.
The Kotschick-Morgan conjecture states that
$$\delta_P(\alpha)=\sum_{i=0}^r a_i(r,d,X) q^{r-i}\alpha^{d-2r-2i}$$
where $q$ is the intersection form of $X$ and the co-efficients $a_i$ depend
only on
$r,d$ and the homotopy type of $X$.  Assuming this conjecture, Gottische,
\cite{G},
has derived formulas for the coefficients $a_i$ in terms of modular forms.

The formulas of Gottische also underlie the relation of the Seiberg-Witten and
Donaldson invariants.  Pidstrigach and Tyurin, in \cite{PT}, propose a program
to compare the two (see \cite{BGP} for a list of references of similar
ideas).  They define a larger moduli space which gives a cobordism between
a link of the anti-self dual moduli space and the link of singularities of
the form $M^{SW}_L\times\Sigma^r(X)$, where $M^{SW}_L$ is a Seiberg-Witten
moduli space.
The link of $M^{SW}_L\times\Sigma^r(X)$ has a description similar to
but more complicated than $D(\alpha,g)$.
Pairing certain cohomology classes with the link of the anti-self dual moduli
space gives a multiple of the Donaldson polynomial.
If we could calculate the pairing of these cohomology classes with the link of
the $M^{SW}_L\times\Sigma^r(X)$, we could give an explicit formula relating the
Seiberg-Witten and Donaldson invariants.
To understand the link of
the singularities $M^{SW}_L\times\Sigma^r(X)$, it would then be helpful to
understand
the \lq\lq easy\rq\rq case of $D(\alpha,g)$.

This paper is organized as follows.  In section one, we give background
information on
an open neighborhood of $[\alpha ]\times\Sigma^r(X)$ and the cohomology
classes.
In section two, we perform the $r=1$ calculation.  The heart of this paper lies
in
section three, the $r=2$ calculation, which is divided into four parts.
Taubes' gluing data description of a neighborhood of $[\alpha
]\times\Sigma^2(X)$
requires two open sets, one for the off-diagonal points and one for the
diagonal points.
In sections 3.1 and 3.2, we give compactifications of these open sets and
compute
pairings.  The pairing with the open set of points on the diagonal, requires
additional
work, done in section 3.3, to describe the moduli space of charge two
instantons on $S^4$.
The final section, 3.4, computes the correction term where the two
compactifications
are compared.

\

\noindent
{\bf 1: Gluing Data Bundles and Links of Reducibles}

\

We consider a family of reducibles, $[\alpha]\times\Sigma^r(X)$, where
$[\alpha]$
is a reducible anti-self dual connection arising from the reduction
$Q_\alpha\times_{S^1}SO(3)$, where $Q_\alpha$ is an $S^1$ bundle with
$c_1(Q_\alpha)=\alpha$
and $r=(\alpha^2-p_1)/4$.
The space $[\alpha]\times\Sigma^r(X)$ is a stratified set, the strata,
$\Sigma$,
of $\Sigma^r X$ are given by partitions of $r$.  If $N=-\alpha^2-2$, then a
neighborhood
of the background connection in the parametrized moduli space
$\Mod {p+4r} \lambda$ is given by $\bC^N/\Gamma_\alpha$.
There is a description due to Taubes, \cite T,
\cite{FM} of a neighborhood of $\bC^N/\Gamma_\alpha\times (\Sigma\backslash
\nu(\Delta))$
where $\nu(\Delta)$ is a neighborhood of the big diagonal (i.e. a neighborhood
of
$\overline\Sigma\backslash \Sigma$ in $\Sigma$).

Let $Z_k(\epsilon)$ be the Uhlenbeck compactification of
the space of centered, charge $k$ instantons
on $S^4$.  That is, if $x:S^4\backslash s\rightarrow \bR^4$ is stereographic
projection from
the south pole,
$$Z_{n_i}(\epsilon)=\{[A]\in {\CMod {n_i} {S^4}}: \int x_i|F_A|^2dx=0,\
\int|x|^2|F_A|^2dx <\epsilon\}.$$
Here, $|F_A|dx$ includes the points of Dirac measure for the generalized
connection $[A]$.
Because the connections in $Z_k(\epsilon)$ are concentrated near the north
pole, we can
extend the space of connections framed at the south pole over this set to get
$Z^s_k(\epsilon)$.
There is an $SO(4)$ action on the space coming from lifting the $SO(4)$
rotation on
$S^4$ to the principal bundle and an $SO(3)$ action on the south pole frame.

If the stratum $\Sigma$ is described by $(n_1,\dots,n_l)$ where $\sum n_i=r$,
then
a neighborhood of $\bC^N/\Gamma_\alpha\times \left( \Sigma\backslash
\nu(\Delta)\right)$ is
described as follows.  Let $Q=Q_\alpha$ be the $S^1$ reduction
corresponding to $\alpha$, $Fr(X)$ be the tangent frame bundle,
$E=Q\times_XFr(X)$, and $\pi_i:X^l\rightarrow X$ projection onto the $i$-th
factor.
$$Gl(\Sigma,\alpha)=\prod_{i=1}^l\pi_i^*E\times_{SO(4)}Z^s_{n_i}(\epsilon)/S^1$$
This is a bundle over $X^l$, to get a bundle over $\Sigma$, we have to mod out
by
$P(n_1,\dots,n_l)$ the subgroup of the symmetric group on $l$ elements
preserving the
multiplicities $n_i$.  Our neighborhood is then
$$\nu_\alpha=\left( \bC^N\times_{\Gamma_\alpha} Gl(\Sigma,\alpha)\right)
/P(n_1,\dots,n_l).$$
$\Gamma_\alpha$ acts by diagonal multiplication on the factors of $Q$ in each
$E$.
This space is fibered by
$\pi_\Sigma: \nu_\alpha\rightarrow\bC^N/\Gamma_\alpha\times (\Sigma\backslash
\nu(\Delta))$.

In \cite{T}, Taubes describes a homeomorphism from these spaces of gluing data,
$\nu_\alpha$
onto their image in the moduli space.
We divide this map into two parts, a splicing map and a pertubation of the
image of
the splicing map onto the moduli space.
The splicing map $\gamma'_\Sigma$ describes how
the background connection in $\bC^N/\Gamma_\alpha$ and the instantons on $S^4$
are
spliced together using cut-off functions and the frames to
give an almost anti-self dual connection.
The perturbation, which we write as $\xi$, moves the image of $\gamma'_\Sigma$
onto the parametrized moduli space.
Note that the element of $\bC^N/\Gamma_\alpha$ describes a connection and a
metric,
so the gluing is performed with respect to the appropriate metric.
The composition, $\gamma_\Sigma=\gamma'_\Sigma+\xi_{\gamma'}$
is called the gluing map.  The union of the images of $\nu_\alpha$
cover a neighborhood of $\bC^N/\Gamma_\alpha\times \Sigma^r(X)$ and we can
define a link of these reducibles in this neighborhood.  This defines a
homology
class $[D(\alpha,N,g)]$ ($g$ being the metric with respect to which $[\alpha]$
is anti-self dual.  The difference term, for $\alpha^2<-1$ is calculated by
$$\delta_P(\alpha)(z)=<\mubar(z),[D(\alpha,N,g)]>.$$
For $\alpha^2=-1$, there is an obstruction to gluing.  We write $e(O)$ for the
Euler
class of this bundle and our calculation can be expressed as
$$\delta_P(\alpha)(z)=<\mubar(z),e(O)\cap [D(\alpha,N,g)]>.$$

On the complement of the reducibles in the image of $\nu_\alpha$,
the $\Gamma_\alpha$ action on
$\bC^N\times Gl(\Sigma,\alpha)$ is free.
We write $c_1$ for the first Chern class of the $\Gamma_\alpha$ action.
For a homology class $x$ on $X$, we write $\Sigma(x)$ to denote
the cohomology class on $\Sigma$ that, under the map $X^n\rightarrow \Sigma$,
pulls back to the symmetrization of $\pi_i^*PD[x]$.
{}From \cite{KoM}, we have the following description of the $\mu$ map
on the images of $\gamma_\Sigma$.
\proclaim{Lemma 1.1} For $x_i\in H_2(X;\bZ)$ and $1\in H_0(X;\bZ)$ the
generator,
$$\align
 \gamma_\Sigma^*\mubar(x) & =
{1\over 2}<\alpha,x>c_1+\pi_\Sigma^*\Sigma (x) \\
\gamma_\Sigma^*\mubar(1) & =
-{1\over 4}c_1^2+\pi_\Sigma^*\Sigma (1) \\
  \gamma_\Sigma^*(e(O)) & =c_1
\endalign$$
\endproclaim

\

The difficulties in computing these pairings arise from two sources.  The first
is the
presence of $Z_k^s(\epsilon)$ which remain mysterious in spite of the ADHM
description.
The second is the problem of the overlap of the sets $\partial \nu_\alpha$.
For $r=0,1$,
neither of these difficulties arise.  The main result here is a method for
dealing with
these problems in the next simplest situation, where two strata are present.

\

\

\noindent
{\bf 2: First Order Reducibles}

\

We begin our work by calculating the case where the reducible has only one
point of
totally concentrated curvature.  This computation has already been performed in
the
$SU(2)$ case in \cite {Y}.  We note that $Z_1^S(\epsilon)\simeq c(SO(3))$.
The $SO(4)=SU(2)_L\times_{\bZ /2}SU(2)_R$ rotation group acts on this via it's
projection
to $SO(3)_R$ (see \cite Y).
We may thus simplify $Fr(X)\times_{SO(4}c(SO(3))\simeq
\Fr\times_{SO(3)}c(SO(3))$
and describe $\nu_\alpha$ as:
$$\nu_\alpha\simeq \Gamma_\alpha\backslash
\left(\bC^N\times Q_\alpha\times_X\Fr\times c(SO(3))/S^1\right)$$
We present a covering of this space by a vector bundle.  Let
$Q_{\alpha\over 2}\rightarrow Q_\alpha$ be a twofold fiberwise covering and
$\Fl\rightarrow \Fr$
be an $SU(2)$ lifting.  Although these bundles need not exist on $X$, one can
pull back
to $X'=X\times_{K(\bZ,2)}K(\bZ,2)$ where $X\rightarrow K(\bZ,2)$ is given by
$\alpha$
and $K(\bZ,2)\rightarrow K(\bZ,2)$ is given by doubling the universal class.
Since
$\alpha$ is an integral lift of $w_2(P)$, pulling back to $X'$ kills $w_2$.
Moreover,
the map $X'\rightarrow X$ is a rational homotopy isomorphism, we can either
pull all
the bundles back to $X'$, work there and divide by the degree, or simply work
on $X$
pretending that these liftings exist.  We thus cover $\nu_\alpha$ with
$$S^1\backslash\left( \bC^N\times Q_{\alpha\over 2}\times_X
\Fl\times_{SU(2)}\bC^2/S^1\right).$$
The degree of this covering map is $-2^N$: the map $\bC^2\rightarrow c(SO(3))$
has degree
$-2$ (see \cite{Oz}), $Q_{\alpha\over 2}\rightarrow Q_\alpha$ has degree $2$
and
to make the covering equivariant with respect to the stabilizer action, we have
to
square each of the $\bC$'s in $\bC^N$.  The stabilizer action doubles,
subtracting one
factor of two from the degree and pulling back $c_1$ of the $\Gamma_\alpha$
action
to $2h$ where $h$ is the $c_1$ of the lefthand $S^1$ action upstairs.  We note
that the
cover of $\partial\nu_\alpha$ can be written as ${\Bbb P}(\bC^N\oplus \tilde
E)$ where
$\tilde E\simeq Q_{\alpha\over 2}\otimes (\Fl\times_{SU(2)}{\overline\bC}^2)$.
If $p_+=p_1(\Fr)$,
then $c_1(\tilde E)=\alpha$ and $c_2(\tilde E)={1\over 4}(\alpha^2-p_+)$.
If we write $b_i=<\alpha,x_i>$, our computation is now
$$\align
&\delta_P(\alpha)(x_1,\dots,x_{d-2l},1^l)
={(-1)\over 2^N}<\biggl(\prod_{i=1}^{d-2l}(b_ih+\pi^*PD[x_i])\biggr)
(-h^2+\pi^*PD[1])^l,
[{\Bbb P}(\bC^N\oplus \tilde E)]>  \tag {2.1} \\
&={(-1)^{l+1}\over
2^N}<(\prod_{j=1}^{d-2l}b_j)h^d+\sum_{i=1}^{d-2l}(\prod_{j\neq
i}b_j)\pi^*PD[x_i]h^{d-1}
+{1\over 2}\sum_{r,s=1}^{d-2l}(\prod_{i\neq r,s}
b_I)\pi^*PD[x_r]\pi^*PD[x_s]h^{d-2} \\
&-l(\prod_{j=1}^{d-2l}b_j)h^{d-2}\pi^*PD[1],
[{\Bbb P}(\bC^N\oplus \tilde E)]>
\endalign$$
We use the definition of characteristic class on this projectivization:
$h^{N+2}-\pi^*c_1(\tilde E)h^{N+1}+\pi^*c_2(\tilde E)h^N=0$ and $d=N+3$ to
rewrite
the above equation as
$$\align
&={(-1)^{l+1}\over 2^N}<(\prod_{j=1}^{d-2l} b_j)(\pi^*\alpha^2-{1\over
4}(\alpha^2-p_+)h^{N+1}
+\sum_{i=1}^{d-2l-1}(\prod_{j\neq i}b_j)\pi^*\alpha\pi^*PD[x_i]h^{N+1} \tag
{2.2}\\
&+{1\over 2}\sum_{r,s=1}^{d-2l}(\prod_{j\neq r,s}b_j)h^{N+1}
\pi^*PD[x_r]\pi^*PD[x_s]
-l(\prod_{j=1}^{d-2l}b_j)h^{N+1}\pi^*PD[1]
,[{\Bbb P}(\bC^N\oplus \tilde E)]>
\endalign$$
Each of these terms contains a factor pulled up from $X$ dual to a point.
Thus, we have
reduced our computation to the pairing $<h^{N+1},{\Bbb CP}^{N+1}>=(-1)^{N+1}$.
We further simplify by noting that $\alpha^2=p_1+4=1-d$. The calculation for
the
case $\alpha^2=-1$ is done similarly.  Our answer is

\proclaim{Proposition 2.1} For a reducible connection $\alpha$ with
$\alpha^2=p_1+4$,
$\alpha^2<-1$
$$\align
\delta_P(\alpha)&(x_1,\dots,x_{d-2l},1^l) \\
&=\left( (-{1\over 2})^{d+l-1}\left((d+p_++3-12l)\alpha^{d-2l}\right)
+{(-1)^{d-1}(d-2l)! \over 2^{d-3}2! (d-2l-2)!} q\alpha^{d-2l-2}\right)
(x_1,\dots,x_{d-2l})\\
&={(-1)^{d+l-1}\over
2^d}\left((2d+2p_++6-24l)\alpha^{d-2l}+4(d-2l)(d-2l-1)q\alpha^{d-2l-2}\right)
(x_1,\dots,x_{d-2l})
\endalign$$
If $\alpha^2=-1$, so $d=2$ we have
$$\align
&\delta_P(\alpha)(x_1,x_2)=\left({1\over 4}\right)\left(
(5+p_+)\alpha^2+4q\right)(x_1,x_2) \\
&\delta_P(\alpha)(1)=-{1\over 4}(2p_+-14)
\endalign $$
\endproclaim

\

{\bf Remark:} This calculation agrees with that of Friedman and Qin in
\cite{FQ, Thm 6.4, 6.5}
for the cases $l=0$ and $l=1$. The co-efficient of $q\alpha^{d-2l-2}$ agrees
with that
of \cite{KoM} up to a sign error in \cite{KoM} as noted in \cite{KoL}.

\

\noindent
{\bf 3: Second Order Reducibles}

\

We now consider the case where $\alpha^2-p_1=8$, where the set of reducibles is
described by $[\alpha]\times \Sigma^2(X)$.  This set is covered by two strata,
the
off-diagonal $\Sigma_{(1,1)}$ and the diagonal $\Sigma_2$.  The off-diagonal,
or
upper stratum, has a neighborhood described by two copies of the bundle
described
in section 2.  We compactify this set by extending the bundles over the
diagonal.
The lower stratum involves the space of charge two instantons on $S^4$,
$Z^s_2(\epsilon)$.
Applying an extension of some work of Hattori's, \cite{Ha,Ha'}, gives enough
information
about the equivariant cohomology of this space to compute the cup-products on
this space.
We give a natural compactification of this set and finally show the error term,
resulting
from the two compactifications, vanishes.  Our result is:

\proclaim{Proposition 3.0} For $\alpha^2=p_1+8$, and $\alpha^2<-1$,
$$\align
(-1)^{d+l}2^d&\delta_P(\alpha)(x_1,\dots,x_{d-2l},1^l) \\
&=\left(2d^2+4dp_++2p_+^2+13d+10p_++21+l(-408-48d-48p_++288l)\right)
\alpha^{d-2l}
(x_1,\dots,x_{d-2l}) \\
&+\left(16p_++16d+32-192l\right)(\alpha^{d-2l-2}q)(x_1,\dots,x_{d-2l}) \\
&+8{(d-2l)!\over (d-2l-4)!}\alpha^{d-2l-4}q^2
\endalign$$
If $\alpha^2=-1$, so that $d=6$ and $0\leq l\leq 3$,
$$\align
(-1)^{l+3}2^6\delta_P(\alpha)&(x_1,\dots,x_{6-2l},1^l) \\
=&\left(2p_+^2-58p_++171-344l+16lp_++160l^2\right) \alpha^{6-2l} \\
&+\left(16p_++128-192l\right)\alpha^{4-2l}q
+8{(6-2l)!\over (2-2l)!}\alpha^{2-2l}q^2
\endalign$$
where the terms with negative exponents will vanish.
\endproclaim
{\bf Remark:} This computation agrees with those of \cite{EG} and \cite{FQ}.

\

{\bf 3.1: The Higher Stratum for r=2}

\

We pull the bundle $Gl(\Sigma_{(1,1)},\alpha)$ back to $X^2$ and extend it over
the diagonal,
writing this compactification as $\nu_{\Sigma_{(1,1)}}\cup C$.
Both factors of $E$ are covered with $\tilde E$.  This branched cover has
degree $2^N$ from
the factors of $\bC$, $2^2$ from the two copies of $E$, $2$ from
$X^2\rightarrow \Sigma^2X$
and ${1\over 2}$ from the doubling of the $\Gamma_\alpha$ action.  Thus, at a
price of a
factor of $2^{-N-2}$ we can pair our cohomology classes with
${\Bbb P}({\Bbb E})={\Bbb P}(\bC^N\oplus\pi_1^*\tilde E\oplus\pi_2^*\tilde E)$
as in the previous section.  In the following computation, we relax our
notation, writing
$x_i$ and $1$ for both the homology classes and their Poincare duals.
$$\align
<\prod_{i=1}^{d-2l}&\mubar(x_i)\cup\mubar(1)^l,[\nu_{\Sigma_{(1,1)}}\cup C ]>
=2^{-N-2}<\biggl(\prod_{i=1}^{d-2l}(b_ih+\pi_1^*x_i+\pi_2^*x_i)\biggr)
(-h^2+\pi_1^*1+\pi_2^*1)^l,
[{\Bbb P}({\Bbb E})]> \tag {3.1.1} \\
&={(-1)^l\over 2^{N+2}}<\biggl(\prod_{j=1}^{d-2l} b_j\biggr)h^d
+h^{d-1}\sum_{i=1}^{d-2l}\biggl(\prod_{j\neq i}
b_j\biggr)(\pi_1^*x_i+\pi_2^*x_i) \\
&+h^{d-2}\left[ \biggl(\prod_{j=1}^{d-2l} b_j\biggr)(-l)(\pi_1^*1+\pi_2^*1)+
\sum_{|I|=2}^d\biggl(\prod_{j\notin I}b_j\biggr)\prod_{i\in
I}(\pi_1^*x_i+\pi_2^*x_i) \right] \\
&+h^{d-3}\left[ \sum_{i=1}^{d-2l}\biggl(\prod_{j\neq i}b_j\biggr)
(\pi_1^*x_i+\pi_2^*x_i)(-l)(\pi_1^*1+\pi_2^*1)+
\sum_{|I|=3} \biggl(\prod_{j\notin I}b_j\biggr)\prod_{i\in
I}(\pi_1^*x_i+\pi_2^*x_i)\right] \\
&+h^{d-4}\left[ \biggl(\prod_{j=1}^{d-2l}b_j\biggr) l(l-1)\pi_1^*1\pi_2^*1
+\sum_{|I|=2}\left(\prod_{j\notin I} b_j\right)\prod_{i\in I}
(\pi_1^*x_i+\pi_2^*x_i)(-l)(\pi_1^*1+\pi_2^*1)\right] \\
& +h^{d-4}\left[\sum_{|I|=4}\biggl(\prod_{j\notin I}b_j\biggr)\prod_{i\in
I}(\pi_1^*x_i+\pi_2^*x_i)\right],
[\bP({\Bbb E})]>
\endalign$$
Here these sums are over all $I\subset \{1,\cdots,d-2l\}$ of the appropriate
size.
Now,
$$\align
c_1({\Bbb E})&=\pi_1^*\alpha+\pi_2^*\alpha \tag {3.1.2}\\
c_2({\Bbb E})&=\pi_1^*\alpha\pi_2^*\alpha+
{1\over 4}(\pi_1^*(\alpha^2-p_+)+\pi_2^*(\alpha^2-p_+)) \\
c_3({\Bbb E})&=\pi_1^*\alpha({1\over 4}\pi_2^*(\alpha^2-p_+))
+\pi_2^*\alpha({1\over 4}\pi_1^*(\alpha^2-p_+)) \\
c_4({\Bbb E})&={1\over 16}\pi_1^*(\alpha^2-p_+)\pi_2^*(\alpha^2-p_+)
\endalign$$
The Leray-Hirsch relationship is, using $p_1=-d-3, \alpha^2=p_1+8,
\alpha^2=-N-2$ to
get $N=d-7$
$$
h^{d-3}-\pi^*c_1({\Bbb E})h^{d-4}+\pi^*c_2({\Bbb E})h^{d-5}-\pi^*c_3({\Bbb
E})h^{d-6}
+\pi^*c_4({\Bbb E})h^{d-7}=0 \tag {3.1.3}$$
This allows us to reduce the terms with $h^k$ in them as in the $r=1$ case,
writing
$c_i$ for $\pi^*c_i(\bC^N\oplus\tilde E\oplus\tilde E)$:
$$\align
&={(-1)^l\over 2^{N+2}}<\biggl(\prod_{j=1}^{d-2l}
b_j\biggr)(c_1^4-3c_1^2c_2+2c_1c_3+c_2^2-c_4)h^{d-4}
+(c_1^3-2c_1c_2+c_3)h^{d-4}\sum_{i=1}^{d-2l}\biggl(\prod_{j\neq i}
b_j\biggr)(\pi_1^*x_i+\pi_2^*x_i)
\tag {3.1.4} \\
&+(c_1^2-c_2)h^{d-4}\left[ \biggl(\prod_{j=1}^{d-2l}
b_j\biggr)(-l)(\pi_1^*1+\pi_2^*1)+
\sum_{|I|=2}^d\biggl(\prod_{j\notin I}b_j\biggr)\prod_{i\in
I}(\pi_1^*x_i+\pi_2^*x_i) \right] \\
&+c_1h^{d-4}\left[ \sum_{i=1}^{d-2l}\biggl(\prod_{j\neq i}b_j\biggr)
(\pi_1^*x_i+\pi_2^*x_i)(-l)(\pi_1^*1+\pi_2^*1)+
\sum_{|I|=3} \biggl(\prod_{j\notin I}b_j\biggr)\prod_{i\in
I}(\pi_1^*x_i+\pi_2^*x_i)\right] \\
&+h^{d-4}\left[ \biggl(\prod_{j=1}^{d-2l}b_j\biggr) l(l-1)\pi_1^*1\pi_2^*1
+\sum_{|I|=2}\biggl(\prod_{j\notin I} b_j\biggr)\prod_{i\in I}
(\pi_1^*x_i+\pi_2^*x_i)(-l)(\pi_1^*1+\pi_2^*1)\right] \\
& +h^{d-4}\left[\sum_{|I|=4}\biggr(\prod_{j\notin I}b_j\biggr)\prod_{i\in
I}(\pi_1^*x_i+\pi_2^*x_i)\right],
[\bP ({\Bbb E})]>
\endalign$$
In the above equation, each term contains an eight dimensional cohomology class
pulled
up from $X^2$.  Cupping with these terms is equivalent to restricting to a
fiber of
${\Bbb P(E)}$ which is ${\Bbb CP}^{N+4}={\Bbb CP}^{d-4}$.  The term $h^{d-4}$
gives
us a factor of $(-1)^{d-4}$, and all that remains is to calculate the
intersection
numbers of the classes pulled up from $X^2$.  We use the notation
$p_+=p_1(\Fr)$.
$$\align
&<\prod_{i=1}^{d-2l}\mubar(x_i)\cup\mubar(1)^l,[\nu_{\Sigma_{(1,1)}}\cup C ]>
\tag{ 3.1.5} \\
&={(-1)^{d+l-4}\over 2^{d-5}}
<\biggl(\prod_{j=1}^{d-2l} b_j\biggr)\left({9\over
16}\pi_1^*\alpha^2\pi_2^*\alpha^2
+{3\over 16}(\pi_1^*\alpha^2\pi_2^*p_+
+\pi_1^*p_+\pi_2^*\alpha^2)+{1\over 16}\pi_1^*p_+\pi_2^*p_+\right) \\
&+\sum_{i=1}^{d-2l} (\prod_{j\neq i} b_j)
\left({3\over 4}(\pi_1^*\alpha\pi_2^*\alpha^2+\pi_1^*\alpha^2\pi_2^*\alpha)
      +{1\over 4}(\pi_1^*\alpha\pi_2^*p_++\pi_1^*p_+\pi_2^*\alpha^2)\right)
(\pi_1^*x_i+\pi_2^*x_i) \\
&+\left[\biggl(\prod_{j=1}^{d-2l}b_j\biggr)(-l)(\pi_1^*1+\pi_2^*1)
+ \sum_{|I|=2}(\prod_{j\notin I}b_j)\prod_{i\in I}(\pi_1^*x_i+\pi_2^*x_i)
\right]
\left[{3\over 4}(\pi_1^*\alpha^2+\pi_2^*\alpha^2)+\pi_1^*\alpha\pi_2^*\alpha
+{1\over 4}(\pi_1^*p_++\pi_2^*p_+)\right] \\
&+\left[\sum_{i=1}^{d-2l}\biggl(\prod_{j\neq i}b_j\biggr)
(\pi_1^*x_i+\pi_2^*x_i)(-l)(\pi_1^*1+\pi_2^*1) +
\sum_{|I|=3}\biggl(\prod_{j\notin I} b_j\biggr)\prod_{i\in
I}(\pi_1^*x_i+\pi_2^*x_i) \right]
\left(\pi_1^*\alpha+\pi_2^*\alpha\right)\\
&+\left[ \biggl(\prod_{j=1}^{d-2l}b_j\biggr) l(l-1)\pi_1^*1\pi_2^*1
+\sum_{|I|=2}\biggl(\prod_{j\notin I} b_j\biggr)\prod_{i\in I}
(\pi_1^*x_i+\pi_2^*x_i)(-l)(\pi_1^*1+\pi_2^*1)\right] \\
&+\left[\sum_{|I|=4}\biggl(\prod_{j\notin I}b_j\biggr)\prod_{i\in I}
(\pi_1^*x_i+\pi_2^*x_i)\right],[X^2]>
\endalign$$
Now, $<\alpha^2,[X]>=p_1+8=5-d$, and using the following combinatorial lemma we
can calculate the intersections on $X^2$.

\proclaim{Lemma 3.1.1} On $X^r$, for $|I|=2r$, $x_i,\alpha\in H_2(X)$,
$$\align
\prod_{i\in I}(\sum_{j=1}^r(\pi_j^*x_i))&=
{(2r)!\over 2^r} q^{(r)}(x_{i_1},\dots,x_{i_{2r}})[X^r] \\
q^r(x_1,\dots,x_{2r-1},\alpha_*)&=(q^{r-1}\alpha)(x_1,\dots,x_{2r-1})
\endalign$$
\endproclaim
When the dust settles, our equation becomes.
$$\align
<\prod_{i=1}^{d-2l}\mubar(x_i)&\cup\mubar(1)^l,[\nu_{\Sigma_{(1,1)}}\cup C ]>
\tag{3.1.6}\\
 &={(-1)^{d+l-4}\over 2^{d-5}}
(\biggl(\prod_{j=1}^{d-2l}b_j\biggr) \biggl[
{9\over 16}(5-d)^2+{3\over 8}(5-d)p_++{1\over 16}p_+^2+{3\over 2}(d-2l)(5-d) \\
&+{1\over 2}(d-2l)p_+-{3\over 2}l(5-d)-{1\over 2}lp_++2{\binom {d-2l} 2}
-2l(d-2l)+l(l-1)\biggr] \\
&+\sum_{|I|=2}\left(\prod_{j\notin I} b_j\right)\left[{3\over
2}(5-d)q(x_{i_1},x_{i_2})
+{1\over 2}q(x_{i_1},x_{i_2})p_+-2lq(x_{i_1},x_{i_2})\right] \\
&+\sum_{|I|=3}\left(\prod_{j\notin I}b_j\right){4!\over
4}q^2(x_{i_1},x_{i_2},x_{i_3},\alpha)
+\sum_{|I|=4}\left(\prod_{j\notin I}b_j\right){4!\over
4}q^2(x_{i_1},x_{i_2},x_{i_3},x_{i_4})
\endalign $$
After this, some work with symmetric polynomials gives us the answer.
Similarly brutal
work gives us a result for reducibles with $\alpha^2=-1$.

\proclaim{Proposition 3.1.2} For a reducible connection with
$p_1+8=\alpha^2<-1$
$$\align
(-1)^{d+l}2^d&<\prod_{i=1}^{d-2l}\mubar(x_i)\cup\mubar(1)^l,[\nu_{\Sigma_{(1,1)}}\cup C ]> \\
&= \left[ 450+28d+60p_++2d^2+4dp_++2p_+^2
+l(-688-48d+-48p_++288l \right] \alpha^{d-2l} (x_1,\dots,x_{d-2l}) \\
& +\left[ 16d+16p_++112-192l\right] {\binom {d-2l} 2}\alpha^{d-2}q
(x_1,\dots,x_d) \\
&+8{(d-2l)!\over (d-2l-4)!}\alpha^{d-2l-4}q^2 (x_1,\dots,x_{d-2l})
\endalign$$
For a reducible with $\alpha^2=-1$, $d=6$
$$\align
(-1)^{l+3}2^6\delta_P&(\alpha)(x_1,\dots,x_{d-2l},1^l) \\
=&\left[ 690-108p_++2p_+^2-624l  +16lp_++160l^2\right]\alpha^{6-2l} \\
&+\left[ 16p_++208-192l\right]{\binom {6-2l} 2}q\alpha^{4-2l}
+{(6-2l)!\over (2-2l)!} 8q^2\alpha^{2-2l}
\endalign$$
where the terms $q\alpha^{4-2l}$ and $q^2\alpha^{2-2l}$ vanish if $l>2$ and
$l>1$,
respectively.
\endproclaim

\

{\bf 3.2: The Lower Stratum}

\

The lower stratum of $[\alpha]\times \Sigma^2X$ has a neighborhood described by
$$\nu_2=\Gamma_\alpha\backslash\left(\bC^N\times \left(Q_\alpha\times_X
Fr(X)\times_{SO(4)}Z_2^s(\epsilon)/S^1\right)\right).$$
On this set, $\mubar(x_i)={b_i\over 2}c_1+2\pi^*x_i$, where $c_1$ is the
first Chern class of the $\Gamma_\alpha$ action, $b_i=<\alpha,x_i>$.
Note that the factor of two results from calculating the restriction to the
diagonal of $1\times\pi_2^*x_i+\pi_1^*x_i\times 1$.  Similarly
$\mubar(1)=-{1\over 4}c_1^2+2\pi^*1$.
Let $T\subset Z_2(\epsilon)$ be the set of generalized connections where
the background connection is trivial.  We also have $T\subset Z^s_2(\epsilon)$
as the $SO(3)$ action on the frame is trivial on this set.  The set of
reducibles in $\nu_2$
is then parametrized by
$$\nu_2(T)=\Gamma_\alpha\backslash \left( Q_\alpha \times_X
Fr(X)\times_{SO(4)}T\right)/S^1
=Fr(X)\times_{SO(4)}T.$$
The set $T=D^4/(\bZ /(2))$ is
not compact, so the link of $\nu_2(T)$ in $\nu_2$ will not be compact.
This non-compactness arises from two singular points moving away from the
diagonal:
i.e. into $\nu_{(1,1)}$ the gluing data associated to the upper stratum.
The link of $\nu_2(T)$ can be presented as the union of two sets.
Let $L^s\subset Z^s_2(\epsilon), L\subset Z_2(\epsilon)$ be the links of $T$ in
these sets.  Our two open sets are then,
$$\align
U_1= & \Gamma_\alpha\backslash\left(S^{2N-1}\times
(Q_\alpha\times_X Fr(X)\times_{SO(4)}Z_2^s(\epsilon)/S^1)\right), \\
U_2= & \Gamma_\alpha\backslash\left( \bC^N\times (Q_\alpha\times_X Fr(X)
\times_{SO(4)}L^s/S^1)\right).
\endalign$$
The set $U_1$ is where the connection on $S^4$ is allowed to approach the cone
point,
but the background connection is kept away from $T$, while the set $U_2$
is where the background connection can approach the reducible, but the $S^4$
connection
may not approach $T$.
The boundary of $T$ is then covered by $U_2$.  We write $\partial L^s$ for the
boundary
of $L^s$.
Now, the $SO(4)$ action on $Z^s_2(\epsilon)$ respects strata and commutes with
the framing
action.  Thus, we can compactify $L^s$ by taking cones on the $SU(2)_L$
orbits in $\partial L^s$.  Write $L^s\cup  C_2^s$, for this compactification.
The framing action extends over this compactification, so $L^s\cup C_2^s$ is an
$SO(3)$ bundle over $L\cup C_2$ the space constructed by taking cones on
the $SU(2)_L$ orbits in $\partial L$.
We then use this compactification to compactify $U_2$,
replacing $L^s$ with $L^s\cup \tilde C_s$.  Write $\overline U_2$ for this
compactification.
We see that the cohomology classes will extend over
this compactification as they can be expressed as characteristic classes of
the equivariant frame bundle $L^s\cup\tilde C_2$.

The pairing of the cohomology classes with $U_1\cup\overline U_2$ can be
written as
$$\align
<\prod_{i=1}^{d-2l}({b_i\over 2}c_1+2\pi^*x_i)&( -{1\over 4}c_1^2+2\pi^*1)^l,
[U_1\cup\overline U_2]>
\tag {3.2.1} \\
=2^{-d+2l}&\bigg[ \bigg(\prod_{j=1}^{d-2l}b_j\biggr)c_1^{d-2l}
	+4\sum_{i=1}^{d-2l}\biggl(\prod_{j\neq i}b_j\biggr)c_1^{d-2l-1}\pi^*x_i \\
&+16\sum_{|I|=2}\biggl(\prod_{j\notin I}
b_j\biggr)c_1^{d-2l-2}q(x_{i_1},x_{i_2})]
(-4)^{-l}(c_1^2-8\pi^*1)^l,[U_1\cup\overline U_2]> \\
={(-1)^l\over 2^d}& <\biggl( \prod_{j=1}^{d-2l}b_j\biggr)c_1^d
+4\sum_{i=1}^{d-2l}\biggl(\prod_{j\neq i}b_j\biggr) c_1^{d-1}x_i \\
&+\biggl(16\sum_{|I|=2}\bigl(\prod_{j\notin I} b_j\bigr)q(x_{i_1},x_{i_2})
-8l\bigl(\prod_{j=1}^{d-2l}b_j\bigr)\biggr)\pi^*1 c_1^{d-2},
[U_1\cup\overline U_2]>
\endalign$$
We now give an argument that we can eliminate the set $U_1$ from our
consideration.
We note that the set $\overline U_2$ can be seen as a vector bundle over
a space $B$:
$$\CD
\bC^N\times_{\Gamma_\alpha}\left(Q_\alpha\times_XFr(X)\times_{SO(4)} (L^s\cup
C_2^s)/S^1\right) \\
@V\pi VV \\
B={\Gamma_\alpha}\backslash \left( Q_\alpha\times_XFr(X)\times_{SO(4)}(L^s\cup
C_2^s)/S^1\right)
\endCD$$
The Thom class of this bundle is $(-c_1)^N$, the negative sign appearing
because
Since $Z^s_2(\epsilon)$ is equivariantly
(with respect to the framing and rotation actions)
retractable onto a point in $T$, $U_1$ retracts onto ${\Bbb CP}^{N-1}\times X$.
Here $c_1$ is the negative of the hyperplane section.
Thus $c_1^N$ vanishes on this open set.
This tells us $-c_1)^N$ is the Thom class of the space $B$.  Formally, we have
\proclaim{Lemma 3.2.1}For a cohomology class $x$ of appropriate dimension,
and $i:B\rightarrow U_1\cup\overline U_2$,
$$<x\cup (-c_1)^N,[U_1\cup\overline U_2]>=<i^*x,[B]>.$$
\endproclaim
Because $d-2>N-1$, we see each term in equation 3.2.1 is divisible by the Thom
class
$(-c_1)^N$ and we can use Lemma 3.2.1.  Now, the
base space $B$
can be rewritten as
$${\Gamma_\alpha}\backslash\left(Q_\alpha\times_XFr(X)\times_{SO(4)}
(L^s\cup C_2^s)/S^1\right)
\simeq Fr(X)\times_{SO(4)}(L^s\cup C_2^s)/S^1$$
Prior to taking the $\Gamma_\alpha$ action, the above could be identified
with the tensor product of $\pi^*Q_\alpha$ with the inverse of the $S^1$ bundle
$Fr(X)\times_{SO(4)}(L^s\cup C_2^s)\rightarrow
Fr(X)\times_{SO(4)}(L^s\cup C_2^s)/S^1$.  We write $h$ for the
first chern class of this latter bundle. The inverse arises
because $\Gamma_\alpha$ acts on the $Q_\alpha$ fibers positively and on this
$S^1$ bundle by inverses via the $S^1$ quotient.
The $c_1$ of the $\Gamma_\alpha$ action then becomes $\pi^*\alpha-h$.
Our pairing has become
$$\align
=&{(-1)^{l+N}\over 2^d} <\biggl(
\prod_{j=1}^{d-2l}b_j\biggr)(\pi^*\alpha-h)^{d-N}
+4\sum_{i=1}^{d-2l}\left(\prod_{j\neq i}b_j\right) (\pi^*\alpha-h)^{d-N-1}x_i
\tag {3.2.2} \\
&+\biggl(16\sum_{|I|=2}\bigl(\prod_{j\notin I} b_j\bigr) q(x_{i_1},x_{i_2})
-8l\bigl(\prod_{j=1}^{d-2l}b_j\bigr)\pi^*1\biggr)(\pi^*\alpha-h)^{d-N-2}
,[Fr(X)\times_{SO(4)}(L^s\cup C_2^s)/S^1]> \\
=&{(-1)^{l+d}\over 2^d}<
\biggl(\prod_{j=1}^{d-2l}b_j\biggr)(h^7-7h^6\pi^*\alpha+21 h^5\pi^*\alpha^2)
-4\sum_{i=1}^{d-2l}\biggl(\prod_{j\neq
i}b_j\biggr)(h^6-6h^5\pi^*\alpha)\pi^*x_i \\
&+\biggl(16\sum_{|I|=2}\bigl(\prod_{j\notin I}b_j\bigr) q(x_{i_1},x_{i_2})
-8l\bigl(\prod_{j=1}^{d-2l}b_j\bigr)\biggr)\pi^*1h^5,
[Fr(X)\times_{SO(4)}(L^s\cup C_2^s)/S^1]>
\endalign$$
We can push these classes down from the $S^1$ quotient to the $SO(3)$ quotient
of
$Fr(X)\times_{SO(4)}(L^s\cup\tilde C_2)\backslash T)$.  This turns $h^r$ into
$2\wp^{r-1\over 2}$ if $r$ is odd and makes $h^r$ vanish if $r$ is even where
$\wp$ is the $p_1$ of the $SO(3)$ action.
Our pairing is now:
$$
\align
&={(-1)^{l+d}\over 2^d}<
\biggl(\prod_{j=1}^{d-2l}b_j\biggr)(2\wp^3+42\wp^2\pi^*\alpha^2)
-4\sum_{i=1}^{d-2l}\biggl(\prod_{j\neq
i}b_j\biggr)(-12\wp^2\pi^*\alpha)\pi^*x_i
\tag {3.2.3} \\
&+\biggl(16\sum_{|I|=2}\bigl(\prod_{j\notin I}b_j\bigl)q(x_{i_1},x_{i_2})
-8l\bigl(\prod_{j=1}^{d-2l}b_j\bigr)\biggr)2\pi^*1\wp^2,
[Fr(X)\times_{SO(4)}(L\cup C_2)]>
\endalign$$

We prove the following proposition in the next section

\proclaim{Proposition 3.3.19} For
$p: ESO(4)\times_{SO(4)}(L\cup C_2)\rightarrow BSO(4)$
and $p_+=p_1+2e$,
$$\align  p_*(\wp^2)&=-{5\over 2} \\
p_*(\wp^3)&=-p_-20p_+
\endalign$$
Pulled back to a manifold with $b^+=1$, $p_*(\wp^3)=48-25p_+$.
\endproclaim
Plugging this in we have
\proclaim{Proposition 3.2.2} If $\alpha^2<-1$,
and for the choice of compactly supported lift of the cohomology class
given by adding on $C_2$,
$$\align
(-1)^{d+l}2^d< &
\prod_{i=1}^{d-2l}\mubar(x_i)\cup\mubar(1)^l,[\partial\nu_{\Sigma_2}]> \\
= & \left((-32p_+-15d-429+280l)\alpha^{d-2l}-80{\binom {d-2l} 2}
\alpha^{d-2}q\right)(x_1,\dots,x_d)
\endalign$$
If $\alpha^2=-1$, then
$$
(-1)^{l+7}2^6\delta_P(\alpha)(x_1,\dots,x_{6-2l},1^l)
=(-50p_+-519+280l)\alpha^{6-2l}-80{\binom {6-2l} 2}\alpha^{4-2l}q
$$
\endproclaim

\

\

{\bf 3.3: ADHM Calculations}

\

In this section we use the ADHM description of $Z^s_2(\epsilon)$
to prove Proposition 3.3.19.
As this section is somewhat involved, we begin with a brief summary.

We decompose the link $L$, of the trivial strata $T$ in $Z_s(\epsilon)$ into
two sets, $A,B$.  The set $A$ is branch covered by $S^3\times\bR^+\times S^4$,
which matches the gluing data bundle associated to $\Sigma_{(1,1)}$ near the
diagonal.  The set $B\simeq (S^3\times S^2)*S^1/S^1$ shows why there is more to
$Z_2(\epsilon)$ than just gluing in two charge one instantons.
In Lemma 3.3.4,
the link is then described as a mapping cone $\partial A\rightarrow B$, with
the factor
of $\bR^+$ in $A$ serving as the cone factor.
Because $\partial L$ is contained in $A$, we are able to describe the cap given
by
taking cones on $SO(3)_L$ orbits from our knowledge of the $SO(4)$ action on
$A$,
the cap turns out to be branch covered by $D^4\times S^4$.
In Lemmas 3.3.5-9, we assemble information about the cohomology of these
spaces.
In Lemmas 3.3.11-15, we calculate the restriction of $\wp$ to $L$ and to the
cap.
Lemma 3.3.16 allows us to take unique lifts of $\wp^2$ to cohomology classes
with compact support on $L$ and on the cap.  The remainder of the section
contains
the calculations using this.

We begin by recalling the ADHM description of the moduli space for charge $k=2$
(see
\cite {A,BoM,Ha}.
$$Z^s_2(\epsilon)=\lbrace\left(\matrix \lambda_1 & \lambda_2\\ a & c\\ c & -a
\endmatrix\right): \lambda_i,a,c\in{\Bbb H},
Im(\overline\lambda_1\lambda_2)=2Im(\overline ca)\rbrace / O(2)$$
Here, $O(2)$ acts by
$$\left(\matrix \Lambda\\ B\endmatrix\right)\rightarrow
\left(\matrix \Lambda T\\ T^tBT\endmatrix\right)$$
We work with the double branched cover given by only modding out by $SO(2)$ at
the expense of a factor of $2$.
We refer to this space as $\tilde Z^s_2(\epsilon)$.
Then the $SO(2)$ action is given by
$$(\matrix \lambda_1 &\lambda_2\endmatrix)\rightarrow
(\matrix\lambda_1 & \lambda_2\endmatrix )
\left(\matrix cos(\theta) & -sin(\theta)\\ sin(\theta) &
cos(\theta)\endmatrix\right)
=(\matrix\lambda_1 & \lambda_2\endmatrix ) R(\theta )$$
and
$$ (\matrix a & c\endmatrix)\rightarrow (\matrix a & c\endmatrix)
\left(\matrix cos(2\theta) & -sin(2\theta)\\ sin(2\theta) &
cos(2\theta)\endmatrix\right)
=(\matrix a & c \endmatrix ) R(2\theta )$$

We note that the ADHM description is usually accompanied by a rank condition.
Relaxing
this rank condition, as done here, is equivalent to compactifying the moduli
space
(see \cite {BoM, Thm. 4.4.3} or \cite{DK, Cor 3.4.10}).
Requiring $tr(B)=0$ is equivalent to the condition
$\int x_i|F_A|^2dx=0$ (see \cite {Ma}).  The ADHM correspondence is equivariant
with respect to the $SO(3)$ action on the frame and the $SO(4)$ rotation action
when
the group actions are given on the ADHM data by (for $q\in SO(3),
(p_L,p_R)\in SU(2)\times_{{\Bbb Z}/ 2} SU(2)_R\simeq SO(4)$)
$$\left(\matrix \lambda_1 & \lambda_2 \\ a & c \\ c & -a\endmatrix \right)
\rightarrow
\left(\matrix q\lambda_1\br p_R & q\lambda_2\br p_R \\
p_La\br p_R & p_L c\br p_R \\ p_L c\br p_R & -p_L a\br p_R \endmatrix \right)$$
To keep our orientation, the standard $SO(4)$ action on the northern hemisphere
of $S^4$, ${\Bbb H}$, is given by $h\rightarrow p_L h \br p_R$.  We will write
$\tilde G$ for $SU(2)_L\times SU(2)_R$.
Note that even on $\tilde Z^s_2(\epsilon)$, the $SO(2)$ action gives an
equivalence
between $((a,c),(\lambda_1,\lambda_2))$ and $((a,c),(-\lambda_1,-\lambda_2)$.
Thus the framing action on $\tilde Z^s_2(\epsilon)$
is only an $SO(3)$ action.  We write $\tilde Z_2(\epsilon)$ for the quotient by
the
$SO(3)$ action on the frame.

The trivial strata $T$ is given by $\{\Lambda=0\}$.  The inclusion
$\{\Lambda=0\}\subseteq T$ follows because $T$ is the set of fixed points of
the
$SO(3)$ action on the frame.  The reverse inclusion follows from the
description in
\cite {BoM, equation 4.4.6} of the lower strata.

\

We are interested in $Z^s_2(\epsilon)/SO(3)$.  A clever trick of Hattori's,
from \cite {Ha} allows us to eliminate the ADHM equation.  He observes that
the $SU(2)$ Hopf fibration map gives us $Im(\br \lambda_1\lambda_2)$.
$$\CD
S^3\backslash {\Bbb H}^2 @>\sim >> {\Bbb C}\times {\Bbb R}^3\\
[\lambda_1,\lambda_2]@>>> \left(Re(\br \lambda_1\lambda_2),
{|\lambda_1|^2-|\lambda_2|^2\over 2}, Im(\br\lambda_1\lambda_2)\right)
\endCD$$
The ADHM equation, $Im(\br\lambda_1\lambda_2)=2Im(\br ca)$ allows us to forget
the ${\Bbb R}^3$ factor.  Better still, this map is $SO(2)$ equivariant, if
$SO(2)$ acts on ${\Bbb C}\times {\Bbb R}^3$ by $(z,x)\rightarrow
(ze^{2i\theta},x)$.
We have shown
\proclaim{Lemma 3.3.1} $\tilde Z_2(\epsilon)\simeq {\Bbb H}^2\times_{S^1}{\Bbb
C}$.
Here $S^1$ acts by $R(\theta)$ on ${\Bbb H}^2$ and by $e^{i \theta}$ on ${\Bbb
C}$.
Under this equivalence, the trivial strata are given by
$T=\{[(a,c),z]: z=0,Im(\br ca)=0\}$.
This homeomorphism is $SO(4)$ equivariant if $SO(4)$ acts on
$\bH^2\times_{S^1}\bC$
by $[h_1,h_2,z]\rightarrow [p_Lh_1\overline p_R,p_Lh_2\overline p_R,z]$.
\endproclaim

\

Analyzing the link of $T$ is easier if we use a further decomposition of this
space.
We apply this same Hopf fibration to the $SU(2)_L$ action and get a map
$$\CD
{\Bbb H}^2\times_{S^1}{\Bbb C} @>\pi>> {\Bbb C}\times {\Bbb R}^3\times_{S^1}
{\Bbb C} \\
[(a,c),z] @>>> \left(Re(\br ca),{|c|^2-|a|^2\over 2},Im(\br ca),z\right)
\endCD$$
In this description, $T=\pi^{-1}([z_1,0,0])$.

\

\proclaim{Lemma 3.3.2} For the map $\pi$ as above and co-ordinates
$(z_1,x,z_2)$
on ${\Bbb C}\times {\Bbb R}^3\times {\Bbb C}$

$T=\pi^{-1}(\{x=z_2=0\})=c(S^3\times S^1)/S^1$

$A=\pi^{-1}(\{z_1\neq 0\})\simeq (S^3\times {\Bbb R}^3\times_{{\Bbb Z}/2}
S^1\times
{\Bbb R}^+)\times_{S^1} {\Bbb C}$ where $S^1$ acts on $S^1$ and on ${\Bbb C}$
by
$[\theta,z]\rightarrow [\theta+\phi,ze^{i\phi}]$.  $\bZ/(2)$ on the factors
$S^3$
and $S^1$ by $(q,\theta)\rightarrow (-q,\theta+\pi)$.

$B=\pi^{-1}(\{z_1=0\})\simeq c(S^3\times S^2)\times_{S^1}{\Bbb C}$ where $S^1$
acts by
$(q,x,z)\rightarrow (qe^{-x\theta},x,ze^{i\theta})$
\endproclaim
\demo{Proof}
For $[(a,c),z]\in T$, such that $\pi([(a,c),z]=(z_1,0,0)$, we have
$\br ca=Re(z_1)$. Thus $(a,c)=q(r,s)$ for $q\in S^3,r,s\in{\Bbb R}$. The
$r,s$ give us the $S^1$ and the cone parameter.

For $[(a,c),z]\in B$, such that $\pi([(a,c),z]=(0,x,z_2)$ we have $\br ca=x$,
a purely imaginary quaternion and thus $x\in \bR^3$
and $|c|=|a|$.  Writing $c=|c|q$ for $q\in S^3$, we then have $a={qx\over |c|}$
so $|a|={|x|\over |c|}$.  Thus $|a|=|c|=\sqrt{|x|}$ and we have
$(a,c)=q({x\over \sqrt{|x|}},\sqrt{|x|})$.  The $SO(2)$ rotation action on
$(a,c)$
translates to
$$\align
& ({qx\over \sqrt{|x|}},q \sqrt{|x|})\rightarrow
q\left( {x\over \sqrt{|x|}}cos(\theta )+\sqrt{|x|}sin(\theta ),
-{x\over \sqrt{|x|}}sin(\theta)+\sqrt{|x|}cos(\theta )\right) \\
& =q\left( {x\over |x|}sin(-\theta )+cos(\theta )\right)
\left({x\over\sqrt{|x|}},\sqrt{|x|}\right)
=qe^{-x\theta}\left( {x\over\sqrt {|x|}},\sqrt{|x|}\right)
\endalign $$
B is then parametrized by, for $q\in S^3,x\in{\Bbb R}^3$
$$\CD
c(S^3\times S^2)\times_{S^1}{\Bbb C}@>\phi_B>>{\Bbb H}^2\times_{S^1}{\Bbb C} \\
[(q,x),z] @>>> [q({x\over \sqrt{|x|}},\sqrt{|x|}),z]
\endCD $$

For $[(a,c),z]\in A$, such that $\pi([(a,c),z])=(z_1,x,z_2)$, we begin by
assuming
$Im(z_1)=0$ so that $|a|=|c|$.  Then $\br ca=Re(z_1)+x$.  Let $c=|c|q$ for
$q\in S^3$
and we have $a={1\over |c|}q(Re(z_1)+x)$.  Then
$|a|={1\over |c|}\sqrt{|Re(z_1)^2+|x|^2}$ so
$|a|=|c|=\root 4\of{|Re(z_1)^2+|x|^2}$. Thus for $Im(z_1)=0$, we have
$$(a,c)=q\left( {Re(z_1)+x\over\root 4\of{|Re(z_1)^2+|x|^2}},
\root 4\of{|Re(z_1)^2+|x|^2}\right)$$
We get rid of the assumption that $Im(z_1)=0$ by using the $SO(2)$ equivariance
of
the Hopf fibration map, that is we rotate.  However, since the rotation on
$(a,c)$
gives $z_1\rightarrow z_1e^{2i\theta}$, we must take an involution into account
before parametrizing $A$ by
$$\CD
S^3\times {\Bbb R}^3\times_{{\Bbb Z}/2} S^1\times {\Bbb R}^+\times_{S^1}{\Bbb
C}@>\phi_A>>
{\Bbb H}^2\times_{S^1}{\Bbb C}
\endCD$$
$$\phi_A([q,x,\theta,r,z])=q\left({\matrix {r+x\over\root 4\of{|r^2+|x|^2}} &
\root 4\of{|r^2+|x|^2}\endmatrix}\right)
\left({\matrix cos(\theta) & -sin(\theta)\\ sin(\theta) & cos(\theta)
\endmatrix}\right).
$$
The ${\Bbb Z}/2$ involution is given by sending $q\rightarrow -q$
and $\theta\rightarrow \theta+\pi$.
\enddemo

\

To complete our description of $Z_2(\epsilon)$, we must see how the
spaces $A$ and $B$ fit together.

\proclaim{Lemma 3.3.3}
$\tilde Z_2(\epsilon)$ can be described as the mapping cone of the map
$$S^3\times {\Bbb R}^3\times_{\bZ/(2)} S^1\times_{S^1}{\Bbb C}\rightarrow
c(S^3\times S^2)\times_{S^1}{\Bbb C}$$
$$m([q,x,\theta,z]=[qe^{{-x\theta\over |x|}},{x\over |x|},|x|,z]$$
where $|x|$ represents the cone parameter in $c(S^3\times S^2)$.
\endproclaim
\demo{Proof}
This follows from
${\underset {r\rightarrow 0}\to {\text { lim
}}}\phi_A([q,x,\theta,r,z])=\phi_B([qe^{-x\theta},x,z])$.
\enddemo
We write $\tilde A$ for the branched covering of $A$ given by omitting the
$\bZ/(2)$
action, so $\tilde A=S^3\times \bR^3\times \bC\times\bR^+$.
Lemma  3.3.2 shows us that the link of $T$ in $A$ is branch covered by
the subset of $\tilde A$ given by
$S^3\times S^4\times \bR^+$ and the link in $B$ is $(S^3\times S^2)*S^1/S^1$.
Putting this together with Lemma 3.3.3, we have the following description of
$L$.
\proclaim{Lemma 3.3.4} The link of $T$ in $Z_2(\epsilon)$, $L$, is degree four
branch covered
by the mapping cone of the  natural map
$$m: S^3\times (S^2*S^1)\rightarrow (S^3\times S^2)*S^1/S^1.$$
The $S^1$ action on $(S^3\times S^2)*S^1$ is given by
$$[q,x,\theta]\rightarrow [qe^{x\phi},x,\theta+\phi].$$
This covering is $SO(4)$ equivariant if $[p_L,p_R]\in SO(4)$ acts by
$$[q,x,z]\rightarrow [p_Lq\br p_R,p_Rx\br p_R,z]$$
where $q\in S^3,x\in \bR^3,z\in\bC$ on B and
on A by, for $\theta\in S^1,r\in\bR^+$):
$$[q,x,\theta,r,z]\rightarrow [p_Lq\br p_R,p_Rx\br p_R,\theta,r,z].$$
\endproclaim
\demo{Proof}
This follows from Lemma 3.3.3.  The equivariance follows by noting that
the parametrizing maps $\phi_A$ and $\phi_B$ are equivariant with
respect to these $SO(4)$ actions.
\enddemo
We write $\tilde L$ for the branch cover of the link $L$ described in Lemma
3.3.4.
This lemma shows us
$\partial L$ is branch covered by $S^3\times S^4\subset \tilde L\cap \tilde A$.
The compactification is $L\cup C$ is given by gluing on cones on the
$SO(3)_L$ orbits.  This action is trivial on $S^4$ and the standard action
on $S^3$.  Thus, we have
\proclaim{Lemma 3.3.5}The compactification $L\cup C$ is degree four branch
covered
by
$$\left( D^4\times S^4\right) \cup
\left( S^3\times S^4\times \bR^+\right)\cup
\left((S^3\times S^2)*S^1/S^1\right).$$
This covering in $SO(4)$ equivariant if the action is as
in Lemma 3.3.4 on the second and third
open sets and $[q,x,z]\rightarrow [p_Lq\overline p_R,p_Rx\overline p_R,z]$ on
the
first open set, for $q\in S^3,x\in S^2,z\in\bC$.
\endproclaim
Our compactification thus decomposes naturally into the mapping cone $\tilde L$
and the cap
$D^4\times S^4$.  The cohomology of the latter is obvious, that of the former
is computed by noting that it retracts onto the image $B\cap L$ and using the
following
lemma.
\proclaim{Lemma 3.3.6} $H^*(B\cap L)\simeq H^*(S^3\times
S^2\times_{S^1}(\bC,S^1))$
for $*>0$.
\endproclaim
\demo{Proof}
The Meyer-Vietoris sequence for the decomposition
$$B\cap L =
\left( S^3\times S^2\times_{S^1}\bC\right) \cup
\left( c(S^3\times S^2)\times_{S^1}S^1\right)$$
is the same as that of the pair $S^3\times S^2\times_{S^1}(\bC,S^1)$
because the second of the two open sets retracts to a point.
\enddemo
Another useful piece of information is the following Euler class computation.
\proclaim{Lemma 3.3.7} (\cite{Ha}) The quotient map of the $S^1$ action on
$S^3\times S^2$ given by
$(q,x)\rightarrow (qe^{-x\theta},x)$ can be written
$(q,x)\rightarrow (qx\br q,x)$. The Chern class of the bundle
$S^3\times S^2\rightarrow S^2\times S^2$ given by $[q,x]\rightarrow [qx\br
q,x]$
is $\eta_L-\eta_R$, where $\eta_L$ is the Poincare dual of $\{p\}\times S^2$
and $\eta_R$ that of $S^2\times \{p\}$ in $S^2\times S^2$.
\endproclaim
\demo{Proof}
If we restrict the bundle to $S^2\times \{i\}$, we have the projection map
$[q,i]\rightarrow [qi\br q,i]$ which is the standard right $S^1$ action on
$S^3$.
On the second $S^2$, that is $\{i\}\times S^2$, the bundle is
$\{[q,x]: qx\br q=i\}=\{[q,x]: \br q iq=x\}$ which is given by a left $S^1$
action
on $S^3$.  Thus the co-efficient of the generator of the right $S^2$ in this
$c_1$
is $-1$.
\enddemo
\proclaim{Definition 3.3.8}
Let $T\in H^2(B\cap L)$ be the class corresponding, under the isomorphism of
Lemma 3.3.6
to the Thom class.
\endproclaim
Thus by Lemma 3.3.7, if $i:S^2\times S^2\rightarrow B\cap L$ is
the inclusion, we have $i^*T=\eta_L-\eta_R$.  To compute the cohomology ring of
the compactification, we need to know $m^*: H^*(B\cap L)\rightarrow
H^*(S^3\times S^4)$.
\proclaim{Lemma 3.3.9} $m^*(T\eta_R)=m^*(T\eta_L)=PD[S^3\times \{p\}]$ for
$p\in S^4$.
\endproclaim
\demo{Proof}
The map $m$ respects the decompositions
$S^3\times S^4=\left( S^3\times S^2\times\bC\right)\cup\left(
S^3\times\bR^3\times S^1\right)$
and
$B\cap L=\left(S^3\times S^2\times_{S^1}\bC\right)\cup\left(c(S^3\times
S^2)\times_{S^1}S^1\right)$.
The four dimensional cohomology of both $S^3\times S^4$ and $B\cap L$ is
generated by
the image of the cohomology of the pairs,
$S^3\times S^2\times (\bC,S^1)$ and $S^3\times S^2\times_{S^1}(\bC,S^1)$.
Thus, it suffices to discuss
$m^*:H^4(S^3\times S^2\times_{S^1}(\bC,S^1))
\rightarrow H^4(S^3\times S^2\times (\bC,S^1))$.
We look at the exact sequences of the pairs, writing $m_0$ for $m$ restricted
to the
boundaries,
$$\CD
S^3\times S^2\times S^1 @> i'>> S^3\times S^2\times \bC @>j'>> S^3\times
S^2\times (\bC,S^1) \\
@V m_0 VV @V m VV @V m VV \\
S^3\times S^2\times_{S^1} @>i>> S^3\times S^2\times_{S^1}\bC @>j>>
S^3\times S^2\times_{S^1} (\bC,S^1)
\endCD$$
Let $\eta'$ be the Poincare dual of
$S^3\times\{i\}\times\bC$ in the upper sequence, for $i\in S^2$, let $T'$
be the dual of $S^3\times S^2\times\{0\}$ in the relative cohomology of the
upper
sequence.
Let $T,\eta_L,\eta_R$ be as definition 3.3.8 and Lemma 3.3.7.
Because $T$ is dual to $S^3\times S^2\times_{S^1}\{0\}$, $m^*T$ is
dual to $S^3\times S^2\times \{0\}$, so $m^*T=T'$.  We see
$m^*(\eta_L-\eta_R)=m^*j^*T=j^*T'=0$.  Finally, we see $T\eta_R$
is dual to $S^3\times \{i\}\times_{S^1}\{0\}$ for $i\in S^2$, so $m^*(T\eta_R)$
is dual to $S^3\times \{i\}\times\{0\}$ so $m^*(T\eta_R)=T'\eta'$.
\enddemo
Now that we have described the cohomology of $\tilde L\cup C$, we
wish to extend this description to an equivariant one.  It seems easier
to work with $\tilde G=SU(2)\times SU(2)$ equivariant cohomology than
with $SO(4)$.  Working with rational coefficients, no information is lost, as
evidenced by
this lemma.
\proclaim{Lemma 3.3.10}\cite{HH} Let $c_L,c_R$ be the second Chern classes of
the bundles
$E\tilde G/SU(2)_R\rightarrow B\tilde G$ and
$E\tilde G/SU(2)_L\rightarrow B\tilde G$ respectively, then
$H^*(B\tilde G;\bQ)=\bQ[c_L,c_R]$.

If $s:B\tilde G\rightarrow BSO(4)$ is the classifying map for the natural
$SO(4)$
bundle over $B\tilde G$ and $e,p_1\in H^4(BSO(4))$ are the universal Euler
class and the
universal Pontrjagin
class, then
$s^*(p_1+2e)=-4c_R$ and $s^*(p_1-2e)=-4c_L.$
\endproclaim
We now describe the equivariant cohomology of the cap $D^4\times S^4$ and
calculate the restriction of the equivariant extension of $\wp$ to this set.
\proclaim{Lemma 3.3.11} If $H\in H^4_{\tilde G}(D^4\times S^4)$ is the
generator of
the cohomology of the fiber of
$$\pi_{\tilde G}: E\tilde G\times_{\tilde G}(D^4\times S^4)\rightarrow B\tilde
G,$$
then
$$H^*_{\tilde G}(D^4\times S^4)\simeq
H^*(B\tilde G)[H]/(H^2+2H\pi^*_{\tilde G}c_R+\pi^*_{\tilde G}c_R).$$
Restricted to $D^4\times S^4$, the equivariant extension of
$\wp$ equals $-4H-4\pi^*_{\tilde G}c_R$.
\endproclaim
\demo{Proof}
The assertion about the structure of the equivariant cohomology ring
follows by noting that $D^4\times S^4$ retracts equivariantly onto $S^4$
and the $\tilde G$ action on $S^4$, $[x,z]\rightarrow [p_Rx\overline p_R,z]$
for $x\in \bR^3,z\in\bC$ is equivalent to the action
$[h_1\overline p_R,h_2\overline p_R]$ on $\bH\bP^1$ ($h_i\in\bH$).
Our homotopy quotient retracts onto the quaternionic projectivization of
$BSU(2)_L\times ESU(2)_R\times_{SU(2)_R}\bH^2$ where the cohomology ring is as
described.

The four dimensional cohomology of this homotopy quotient, then has three
generators
$H,c_R,c_L$.  We calculate the co-efficient of $H$ in $\wp$ by restricting to a
single
fiber.  We calculate the co-efficients of $c_R,c_L$ by restricting to the
section of
$\pi_{\tilde G}$ given by $E\tilde G\times_{\tilde G}[1,0,e^0]$, for
$1\in S^3, 0\in\bR^3,e^0\in\bC$.

Restricted to a single fiber, we calculate the bundle as follows.
The Hopf fibration in Lemma 3.3.1, and thus the fibration associated to the
$SO(3)$ framing
action, can be seen as the restriction of
$\bH^2\times_{S^1}\bH^2\rightarrow \bH^2\times_{S^1} B^5$ to the image of
$$\CD
\bH^2\times_{S^1}\bC @>>> \bH^2\times_{S^1}B^5 \\
[(a,c),z]@>>> [(a,c),({1\over 2}Im(\br c a),z)]
\endCD$$
Of course, the Hopf fibration is not a fibration over $0\in B^5$, but the
intersection
of the image of $\bH^2\times_{S^1}\bC$ with $\bH^2\times \{0\}$ is simply $T$.
Composing the above inclusion with $\phi_A$
we get an inclusion $S^3\times S^4\rightarrow \bH^2\times_{S^1}B^5$ by
$$\CD
S^3\times B^5 @>>> \bH^2\times_{S^1}B^5 \\
[q,(x,z)] @>>> \left[q\left( {1+x\over {\root 4\of{1+|x|^2}}},{\root 4\of
{|1+|x|^2}}\right),
(x,z)\right]
\endCD$$
Here $x\in \bR^3,z\in \bC$.
The Hopf fibration then pulls back to the Hopf fibration on the factor $S^4$,
so the associated $SO(3)$ bundle has $p_1$ equal to $-4$ times the generator of
$H^4(S^4)$.

Then, when we restrict to the section, $E\tilde G\times_{\tilde G}(1,0,e^0)$,
we see $(\lambda_1,\lambda_2)$ must satisfy
$\overline\lambda_1\lambda_2=1, |\lambda_1|=|\lambda_2|$ and can thus be
written
$(\lambda_1,\lambda_2)=(q,q)$ for some $q\in S^3$. The bundle is then $SO(3)$
bundle associated to
$$E\tilde G\times_{\tilde G} S^3\rightarrow B\tilde G$$
where $\tilde G$ acts on $S^3$ by $q\rightarrow q\overline p_R$.
This has $p_1=-4c_R$.
\enddemo

The mapping cone of Lemma 3.3.4 retracts equivariantly onto its image
and we can repeat the computations of Lemma 3.3.6, equivariantly.
\proclaim{Lemma 3.3.12} If $\eta_L,\eta_R\in H^2_{\tilde G}(S^3\times S^2/S^1)$
are the generators of the cohomology of the fiber
of $E\tilde G\times_{\tilde G}(S^3\times S^2)/S^1\rightarrow B\tilde G$
where $\tilde G$ acts on $S^3\times S^2$ and $(S^3\times S^2)/S^1$
as in Lemma 3.3.4, then
$$\align
H^*_{\tilde G}(S^3\times S^2/S^1) & \simeq
H^*(B\tilde G)[\eta_L,\eta_R]/(\eta_L^2+\pi^*c_L,\eta_R^2+\pi^*c_R)\simeq
\bQ[\eta_L,\eta_R], \\
H^*_{\tilde G}(S^3\times S^2) & \simeq
H^*(B\tilde
G)[\eta_L,\eta_R]/(\eta_L-\eta_R,\eta_L^2+\pi^*c_L,\eta_R^2+\pi^*c_R) \simeq
\bQ[y]
\endalign$$
where $y$ is a two dimensional cohomology class.
$$H^*_{\tilde G}(B\cap L)\simeq \text{Ker}
\left[H^*_{\tilde G}(S^3\times S^2/S^1)\oplus H^*(B\tilde G)
\rightarrow
H^*_{\tilde G}(S^3\times S^2)\right].$$
\endproclaim
\demo{Proof}
The first assertion follows by noting that the projection
$$S^3\times S^2\rightarrow S^3\times S^2/S^1$$
where $S^1$ acts by $[qe^{-x\theta},x]$
is given by $(q,x)\rightarrow (qx\overline q,x)$.
Thus the $\tilde G$ action on $S^3\times S^2/S^1=S^2\times S^2$ is given by
$(x_1,_2)\rightarrow (p_L x\overline p_L,p_R x\overline p_R)$.
The homotopy quotient of this action on $S^2\times S^2$
is the product of two copies of $ESU(2)/S^1$ which gives the first assertion.

The second assertion follows from noting that the homotopy quotient of
$S^3\times S^2$ is an $S^1$ line bundle over that of
$S^3\times S^2/S^1$ and the Euler class of this bundle is given by Lemma 3.3.7.
Note that $H^2_{\tilde G}(S^2\times S^2)$ is generated on a single fiber, so
to calculate the Euler class of a line bundle, it suffices to restrict our
attention to a single fiber and the non-equivariant computation holds for
the homotopy quotient.

The third assertion follows from the Meyer-Vietoris sequence and noting that
the cohomology of the intersection, $H^*_{\tilde G}(S^3\times S^2)$ vanishes in
odd dimensions to give injectivity.
\enddemo
By virtue of this injectivity, to define global cohomology classes on the
$\tilde L$, it suffices to give their images in
homotopy quotients of $S^3\times S^2\times_{S^1}\bC$ and
$c(S^3\times S^2)\times_{S^1}S^1$.
\proclaim{Definition 3.3.13} The equivariant cohomology classes $T,K_L,K_R$
are defined by
$$\align
T \rightarrow & (\eta_L-\eta_R,0) \in H^2_{\tilde G}(S^3\times S^2/S^1)\oplus
H^2(B\tilde G)\\
K_L\rightarrow & (c_L-\eta_L\eta_R,0)
\in H^4_{\tilde G}(S^3\times S^2/S^1)\oplus H^4(B\tilde G) \\
K_R\rightarrow & (\eta_L\eta_R-c_R,0)
\in H^4_{\tilde G}(S^3\times S^2/S^1)\oplus H^4(B\tilde G)
\endalign$$
\endproclaim
A note on the origins of $K_L$: fiberwise, this class comes from $T\eta_L$.  We
compute
some useful relations among these classes.
\proclaim{Lemma 3.3.14}
$$\align
T^2 & \rightarrow (-c_L-2\eta_R\eta_L-c_R,0) \\
K_L^2 &=-T^2c_L \qquad K_R^2 =-T^2c_R \\
\endalign$$
\endproclaim
\demo{Proof}
Just observe these are equal when restricted to both open sets.
\enddemo

\proclaim{Lemma 3.3.15}(\cite{Ha}) Restricted to $\tilde L$,
$\wp=T^2-4K_R-4c_R$.
\endproclaim
\demo{Proof}
To compute, we restrict to the two open sets given by the homotopy quotients
of $m^{-1}(S^3\times S^2\times_{S^1}\bC)$ and
$m^{-1}(c(S^3\times S^2)\times_{S^1}S^1)$.

The set $S^2\times S^2\subset L\cap B$, can be written as
$[q,x,0], q\in S^3,x\in S^2$.  The parametrizing map $\phi_B$
takes these points to $[(qx,q),0]\in \bH^2\times_{S^1}\bC$.
Over this set, $(\lambda_1,\lambda_2)$
satisfy $\br\lambda_1\lambda_2=x$ and $|\lambda_1|=|\lambda_2|$. We have
$(\lambda_1,\lambda_2)=q'(1,x)$.  The $SO(2)$ action on these $\lambda_i$ is
given by $q'(1,x)\rightarrow q'e^{x\theta}(1,x)$. To get the $SO(3)$ bundle,
we map $q'\rightarrow Ad(q')$. The  frame bundle over
$S^2\times S^2$ is then
$$(S^3\times S^2)\times_{S^1}SO(3)\rightarrow (S^3\times S^2)/S^1$$
where $S^1$ acts by $[q_1,x,Ad(q_2)]\rightarrow
[q_1e^{-2x\theta},x,Ad(q_2e^{x\theta}]$.
If we denote the $S^1$ bundle
$S^3\times S^2$ over $S^2\times S^2$ given in Lemma 3.3.6 by $V$, then Hattori
shows that our $SO(3)$ bundle is the $SO(3)$ extension of the $S^1$ bundle
given by tensoring
$V^{-1}$ with the $S^1$ bundle given by pulling back $SO(3)\rightarrow S^2$
($c_1=2$) from the right $S^2$.  The argument for this goes as follows.
This tensor product of $S^1$ bundles, $V^{-1}\otimes\pi_R^*SO(3)$, can
be described as the quotient of $S^3\times SO(3)$ by the $S^1$ action
$[q_1,Ad(q_2)]\sim [q_1e^{2q_2i\br q_2\theta},Ad(q_2e^{i\theta})]$.
The projection of this space to $S^2\times S^2$ is given by
$[q_1,Ad(q_2)]\rightarrow [q_1q_2i\br q_2\br q_1,q_2i\br q_2]$.
The $S^1$ action on this $S^1$ bundle is then given by
$[q_1,Ad(q_2)]\rightarrow
[q_1e^{-x\theta},Ad(q_2)]\sim[q_1,Ad(q_2)e^{i\theta})]$.
The $SO(3)$ extension of $V^{-1}\otimes \pi_R^*SO(3)$ is given by
$[q_1,Ad(q_2),Ad(q_3)]\sim[q_1e^{-2x\theta},Ad(q_2),Ad(q_3e^{i\theta})]
\sim[q_1,Ad(q_2e^{i\theta}),Ad(q_3e^{i\theta})]$.
Finally, Hattori gives the
map $V^{-1}\otimes\pi_R^*SO(3)\times SO(3)/S^1\rightarrow S^3\times S^2\times
SO(3)/S^1$
by: $[q_1,Ad(q_2),Ad(q_3)]\rightarrow [q_1,q_2i\br q_2,Ad(q_3\br q_2)]$.
One checks that this map respects the $S^1$ action and
$\wp$ and $(\eta_L-3\eta_R)^2$ are equal on the first open set.
We see $(T^2-4K_R-4c_R)$ also maps to $(\eta_L-3\eta_R)^2$ when restricted this
set
(use $T-2\eta_R=\eta_L-3\eta_R$).

The fiberwise contractible open set $E\tilde G\times_{\tilde G}c(S^3\times
S^2)\times_{S^1}$
is easier.  The space $c(S^3\times S^2)\times_{S^1}S^1$
retracts $\tilde G$ equivariantly onto the cone point.
Taking the point in $S^1$ to be $1$, so $\overline\lambda_1\lambda_2=1$ and
$\lambda_1=\lambda_2\in S^3$,
the bundle over the cone point is simply
$E\tilde G\times_{\tilde G}S^3$, where $\tilde G$ acts on $S^3$ by
$q\rightarrow p_Rq$.  Again $T-4K_R-4c_R$ has the same image under restriction
to this open set.
\enddemo
We have calculated the restriction of $\wp$ to two open sets
covering the compactification.  We now observe that this
suffices to give the global cohomology class.  In general,
the possible global extensions of these two restrictions
would be an affine space, with underlying vector space given
by the image of the cohomology of the intersection
of the two open sets under the boundary map of the Meyer-Vietoris
sequence.  We see, however, that our insistence upon working
equivariantly eliminates this ambiguity.
\proclaim{Lemma 3.3.16} $H^*_{\tilde G}(S^3\times S^4)$ vanishes
in odd dimensions.
\endproclaim
\demo{Proof}
We use the Meyer-Vietoris sequence for the decomposition
$$E\tilde G\times_{\tilde G}(S^3\times S^4)=
\left( E\tilde G\times_{\tilde G}(S^3\times S^2\times \bC\right)
\cup
\left( E\tilde G\times_{\tilde G}(S^3\times \bR^3\times S^1\right).$$
The cohomology ring $H^*_{\tilde G}(S^3\times S^2\times\bC)$ has
been calculated in Lemma 3.3.9 and can be seen to be a polynomial
ring with one two dimensional generator, $\bQ[y]$.
The second open set retracts onto $(E\tilde G\times_{\tilde G} S^3)\times S^1$
which has cohomology ring $H^*(B\tilde G)/(c_L+c_R))\otimes
\bQ[\theta]/(\theta^2)$,
where $\theta$ is one dimensional.
The intersection of the two open sets is
$E\tilde G\times_{\tilde G}(S^3\times S^2)\times S^1$
which has cohomology ring $\bQ[y']\otimes \bQ[\theta']/((\theta')^2)$
where again $\theta'$ is one dimensional.
We see that the one dimensional
cohomology class $\theta$ maps onto the one dimensional cohomology of the
intersection,
$\theta'$, and $\bQ[y]$ maps onto $\bQ[y']$,
so none of the odd dimensional cohomology is global.
\enddemo
Thus, our restriction of the equivariant extension of $\wp$ to the two open
sets
$E\tilde G\times_{\tilde G}(D^4\times S^4)$ and
$E\tilde G\times_{\tilde G}\tilde L$ have unique lifts to cohomology classes
with
compact support.
We define the relative pushforward maps by
$$\align
p_1 & : E\tilde G\times_{\tilde G}\left( (D^4,S^3)\times S^4\right)
\rightarrow B\tilde G \\
p_2 & : E\tilde G\times_{\tilde G}\left( B\cap L, S^3\times S^4\right)
\rightarrow B\tilde G.
\endalign$$
We then calculate
$$p_*(\wp^r)=(p_1)_*\left(
(-4H-4c_R)^r\right)+(p_2)_*\left(T^2-4K_R-4c_R)^r\right).$$
We note that, because the dimension of the fundamental class of the fiber is
$8$,
we can restrict to $E\tilde G_{r-8}=\pi_{\tilde G}^{-1}B\tilde B_{r-8}$,
where $B\tilde G_{r-8}$ is the $r-8$ skeleton of some CW decomposition of
$B\tilde G$.
We see the first term will vanish. We use this restriction property
and the relation $H^2+c_RH+c_R^2=0$ to calculate.
$$\align
(p_1)_* & ( (-4)^2(H^2+2Hc_R+c_R^2))=16(p_1)_*(H^2)=16(p_1)_*(-2c_RH+c_R^2)=0
\\
(p_1)_* &((-4)^3(H^3+3H^2c_R+3Hc_R^2+c_R^3))
=-64(p_1)_*\left( H(-2c_RH-c_R^2)+3c_R(-2c_RH+c_R^2)\right)=0.
\endalign$$
To compute the second term, involving $(p_2)_*$, we need the following lemma.
\proclaim{Lemma 3.3.17}
Let $i$ be the inclusion of  $E\tilde G\times_{\tilde G}(S^2\times S^2)$
into the homotopy quotient of the mapping cone of Lemma.
Let $p_3:E\tilde G\times_{\tilde G}(S^2\times S^2)\rightarrow B\tilde G$
be the projection. Then for any
$x\in H^*_{\tilde G}(B\cap L,S^3\times S^4)$, we have
$$(p_2)_*(T^2x)=(p_3)_*(i^*x).$$
The pushforward map $(p_3)_*$ is division by $\eta_L\eta_R$.
\endproclaim
\demo{Proof}
This is equivalent to the assertion that $T^2$ is the Poincare dual of the
space
$E\tilde G_{r}\times_{\tilde G}(S^2\times S^2)$ for any finite $r$.
A neighborhood of this space in the homotopy quotient of the mapping cone is
given by
$$E\tilde G_{r}\times_{\tilde G}(S^3\times S^2\times_{S^1}\bC^2)$$
where one factor of $\bC$ comes from the inclusion of $S^2\times S^2$ into
$B\cap L$
and the other factor arises from the mapping cone.  One can see this
by noting that restricted to the open
set of $B\cap L$ given by $S^3\times S^2\times_{S^1}\bC$, the map $m$ is given
by
$$m: S^3\times S^2\times S^1\times_{S^1}\bC\rightarrow S^3\times
S^2\times_{S^1}\bC.$$
The mapping cone of this restriction of $m$ is then $S^3\times
S^2\times_{S^1}\bC^2$.
The class $T^2$ restricts to this open set of the mapping cone to be the Thom
class
of the normal bundle of $E\tilde G_{r}\times_{\tilde G}(S^2\times S^2)$ and
vanishes on the complement.
\enddemo
We are now ready to perform our computation.  We use
the restriction property: for calculating $(p_2)_*(\wp^r)$ we may
restrict to $E\tilde G_{4r-8}\times_{\tilde G}\tilde L$ so any product
$c_R^ic_L^j$ will vanish if $i+j\geq r-2$.
We also recall from the definition of $T,K_R$ that $i^*T=\eta_L-\eta_R$ and
$i^*K_R=\eta_L\eta_R+c_R$, where $i$ is as in Lemma 3.3.17.
Finally, we employ the cohomology relations of Lemma 3.3.14.
We compute:
$$\align
(p_2)_* & \left( (T^2-4K_R-4c_R)^2\right)=(p_2)_*(T^4-8T^2K_R+16K_R^2)
=(p_2)_*\left( T^2(T^2-8K_R-16c_R)\right)\\
& =(p_2)_*\left( T^2(T^2-8K_R)\right)
=(p_3)_*(i^*(T^2-8K_R))=(p_3)_*(-c_L-2\eta_L\eta_R-c_R-8\eta_L\eta_R-8c_R) \\
&=-10 \\
(p_2)_* & \left( (T^2-4K_R-4c_R)^3\right)
=(p_2)_*\left( (T^2-4K_R)^3-12(T^2-4K_R)^2c_R\right) \\
 & =(p_2)_*\left( T^6-12T^4K_R+48T^2K_R^2-64K_R^3
-12T^4c_R+96T^2K_Rc_R-12(16)K_R^2c_R\right) \\
& =(p_2)_*\left(T^2(
T^4-12T^2K_R+48K_R^2+64K_Rc_R-12T^2c_R+96K_Rc_R-12(16)c_R^2)\right) \\
& =(p_3)_*i^*((\eta_L-\eta_R)^4-12(\eta_L-\eta_R)^2(\eta_L\eta_R+c_R)
+48(\eta_L\eta_R+c_R)^2+64(\eta_L\eta_R+c_R)c_R \\
& \qquad -12(\eta_L-\eta_R)^2c_R+96(\eta_L\eta_R+c_R)c_R) \\
&= (p_3)_*( c_L^4+4c_L\eta_L\eta_R+6c_Rc_L+4\eta_L\eta_Rc_R+c_R^2
-12(-c_L\eta_L\eta_R-3c_R\eta_L\eta_R) \\
& \qquad +48(c_Lc_R+2\eta_L\eta_Rc_R+c_R^2)
 +64(\eta_L\eta_Rc_R+c_R^2) \\
& \qquad -12(-c_Lc_R-2\eta_L\eta_Rc_R-c_R^2)
+96(\eta_L\eta_Rc_R+c_R^2)) \\
& = (4+12)c_L+(4+36+96+64+24+96)c_R =16c_L+320c_R
\endalign$$

Dividing by the degree of the branched covering, $4$, have performed the
following
computation.
\proclaim{Lemma 3.3.18}
$$
p_*(\wp^2)  =-{5\over 2},\qquad
p_*(\wp^3)  = -s^*(p_1-2e)-20s^*(p_1+2e).$$
\endproclaim
Pulling back to $X$, where $b^1(X)=0,b^+(X)=1$, we can see
$$\align
p_+(X) &=p_1(X)+2e(X)=3(b^+-b^-)+4+2b^++2b^-=9-b^-, \\
p_-(X) &= p_1(X)-2e(X)=3(b^+-b^-)-4-2b^+-2b^-=-3-5b^-.
\endalign$$
{}From this, we can derive $p_-(X)=-48+5p_+(X)$
This gives us
\proclaim{Proposition 3.3.19} For the fibration
$p: Fr(X)\times_{SO(4)}(L\cup C_2)\rightarrow X$
where $X$ satisfies $b^1(X)=0, b^+(X)=1$,
if we denote $p_+=p_1(X)+2e(X)$, we have
$$p_*(\wp)^2)=-{5\over 2} \qquad p_*(\wp^3)=48-25p_+$$
\endproclaim

\

{\bf 3.4: The Error Term}

\

We covered a neighborhood of the family of reducibles
$[\alpha]\times \Sigma^2(X)$ with two open sets
$\nu_{\Sigma_{(1,1)}}$ and $\nu_{\Sigma_2}$.  We compactified
$\nu_{\Sigma_{(1,1)}}$ by extending the gluing data bundles over the
diagonal, creating a new space $\nu_{\Sigma_{(1,1)}}\cup C_1$, and
compactified $\nu_{\Sigma_2}$ by gluing on cones on $SO(3)_L$ orbits,
creating the space $\nu_{\Sigma_2}\cup C_2$.
Although the caps do not represent actual connections, we persist in using
the term reducible to describe those points in the compactified spaces of
gluing data where the $\Gamma_\alpha$ action is not free.
We write $(\nu_{\Sigma_{(1,1)}}\cup C_1)^\circ$ and
$(\nu_{\Sigma_2}\cup C_2)^\circ$ for the links of the reducibles in these
space.
The computation can then be expressed as:
$$<\mubar(z),[\partial\nu_\alpha]>=<\mubar(z),[(\nu_{\Sigma_2}\cup C_2)^\circ]>
+<\mubar(z),[(\nu_{\Sigma_2}\cup C_2)^\circ]>
-<\mubar(z),[(C_1\cup C_2)^\circ]>.
\tag {3.4.1}$$
The first two terms on the right hand side were computed in sections 3.1 and
3.2.
Our objective in this section is to show that the final term of this equation
vanishes.

\

The cap $C_1$ is obtained by extending the gluing data bundles over the
diagonal
$\Delta\subset\Sigma^2(X)$.  We then have a natural homeomorphism:
$$k: \bC^N\times_{\Gamma_\alpha}
\left( Q_\alpha\times_X  Fr(X)\times_{SO(4)}
\left( D^4\times c(SO(3))\times c(SO(3))\right)/\bZ/(2)\right)/S^1
\rightarrow C_1.$$
Here the $S^1$ acts on $Q_\alpha$ and $c(SO(3))\times c(SO(3))$ while
$\bZ/(2)$ acts by $-1$ on $D^4$ and by permuting the factors of $c(SO(3))$.
The $SO(4)$ action is the standard action on $D^4$ and $SO(3)_R$ on each factor
of
$c(SO(3))$.
The homeomorphism $k$ is given by noting that a neighborhood of
$\Delta\subset\Sigma^2(X)$
is described by $Fr(X)\times_{SO(4)}D^4/(\bZ/(2))$.  The frames of $Q_\alpha$
and
$Fr(X)$ are parallel translated to frames at the point in $\Sigma^2(X)$
represented by the point of $Fr(X)\times_{SO(4)}D^4/(\bZ/(2)$.
In the case $X=\bR^4$, we could think of this as parallel translating the
frames to the points
$y$ and $-y$, where $y$ is the point in $Fr(X)\times_{SO(4)}D^4/(\bZ/(2)$.
The product
$c(SO(3))\times c(SO(3))$ represents two charge one framed instantons on $S^4$.
We write $C_1^s(x)$ for the cap on a fiber:
$\left( D^4 \times c(SO(3))\times c(SO(3))\right)/\bZ/(2)$.
The reducibles in this cap are given by
$$ \{0\}\times_{\Gamma_\alpha}\left(Q_\alpha\times_X Fr(X)\times_{SO(4)}
\left( D^4\times c\times c\right)/\bZ/(2)\right) /S^1
$$
where $c\in c(SO(3))$ is the cone point.

\

The cap $C_2$ is described by attaching cones on $SO(3)_L$ orbits to
$\partial Z^s_2(\epsilon)$.
We write $C_2^s(x)$ for the space obtained by these cones.
Because the rotation action on $Z^s_2(\epsilon)$ respects the stratification,
the stratification extends over this cap.
The cap $C_2$ is then given by
$$\bC^N\times_{\Gamma_\alpha}\left(Q_\alpha\times_X
Fr(X)\times_{SO(4)}C_2^s(x)\right)/S^1$$
If we write $\overline T$ for
the cones on $SO(3)_L$ orbits in $Z^s_2(\epsilon)$ where the background
connection is trivial,
the reducibles in this cap will be given by
$$\{0\}\times_{\Gamma_\alpha}\left(Q_\alpha\times_X Fr(X)\times_{SO(4)}
\overline T\right)/S^1.$$

We distinguish between the splicing map, where trivializations and cut-off
functions are used
to splice connections on $S^4$ to the background connections on
$Q_\alpha\times_{S^1}SO(3)$,
and the gluing map, where a perturbation is added to the spliced connection to
make it
anti-self dual.  We write $\gamma_{(1,1)}'$ for the splicing map associated to
the gluing
data bundle of the off-diagonal and $\gamma_{(1,1)}$ for the gluing map.
Similarly,
we write $\gamma_2'$ for the splicing map of a charge two $S^4$ instanton along
the diagonal
and $\gamma_2$ for the gluing map.
Then $\gamma_2^{-1}\circ\gamma_{(1,1)}(\nu_{\Sigma_{(1,1)}})\subset
\nu_{\Sigma_2}$
has a boundary with two components, one from the boundary of $\nu_{\Sigma_2}$
and
one from the image of the boundary of $\nu_{\Sigma_{(1,1)}}$.  Because
$\nu_{\Sigma_{(1,1)}}$ and
$\nu_{\Sigma_2}$ give an open cover, these two boundaries are disjoint.  If we
attach $C_1$ and $C_2$ to the appropriate boundaries, we obtain a space we call
$C_2\cup_{\gamma_2^{-1}\circ\gamma_{(1,1)}}C_1$.  The gluing maps respect
the reducibles, so it makes sense to write the link of the reducibles in this
space
as $(C_2\cup_{\gamma_2^{-1}\circ\gamma_{(1,1)}}C_1)^\circ$.
Our error term is then
$$< \mubar(z),[(C_2\cup_{\gamma_2^{-1}\circ\gamma_{(1,1)}} C_1)^\circ]>.$$
This transition map $\gamma_2^{-1}\circ\gamma_{(1,1)}$ is quite mysterious.
We show here that it can be deformed to a transition map $\rho$ which, while
still
not fully understood, has the property that it respects the fibration of the
caps
down to $\bC^N/\Gamma_\alpha\times X$.  With only this information and our work
on
the equivariant cohomology of $\partial L$, we are able to show the error term
vanishes.

\proclaim{Proposition 3.4.1} Let $\rho_x$ be a family of embeddings
parametrized by $x\in X$
$$\rho_x: \left( S^3\times c(SO(3))\times c(SO(3))\right) /\bZ/(2)
\rightarrow Z^s_2(\epsilon),$$
which are $SO(4)$ and $SO(3)$ equivariant with respect to the
rotation and framing actions on the image and the action
$$[q,(A_1,t_1),(A_2,t_2)]\rightarrow
[p_L q\overline p_R, (rA_1\overline p_R,t_1),(rA_2\overline p_R,t_2)],$$
on $\left( S^3\times c(SO(3))\times c(SO(3))\right)/\bZ/(2) $
for $[p_L,p_R]\in SO(4), r\in SO(3),q\in S^3,[A_i,t_i]\in c(SO(3))$.
If we also assume that $\rho_x$ sends
$S^3\times c\times c$ to the trivial strata in $Z^s_2(\epsilon)$
and we glue together $C_1$ and $C_2$ in a fiberwise manner to get the space
$$C_2\cup_\rho C_1 =\bC^N\times_{\Gamma_\alpha}\left( Q_\alpha\times_X
Fr(X)\times_{SO(4)}
(C_2^s(x)\cup_{\rho_x}  C_1^s(x))\right)/S^1,$$
then we can define a link of the reducibles in $C_2\cup_\rho C_1$, which we
write
as $(C_1\cup_\rho C_2)^\circ$.
We can
then calculate that $$<\mubar(z),[(C_2\cup_\rho C_1)^\circ]>=0.$$
\endproclaim
\demo{Proof}
This follows from Lemma 3.3.14 and the computation of $(p_1)_*(\wp^r)$ in
Section
3.3.  Arguing as in Section 3.2, we reduce the computation to an equivariant
computation on $ESO(4)\times_{SO(4)}(C_2^s(x)\cup_{\rho_x}
C_1^s(x))^\circ)/SO(3)$
of pushforwards of $\wp^r$, where now, $\wp$ is Pontrjagin class of the $SO(3)$
action here.  We can give branched covers of $C^s_i(x)/SO(3)$ by $D^4\times
S^4$.
Again, the odd dimensional equivariant cohomology of the intersection
of the two caps vanishes so we have unique lifts of $\wp$ to fiberwise
compactly
supported equivariant cohomology and the exact same computations give our
result.
\enddemo

In the remainder of this section, we describe such a transition map $\rho$ and
then show that $\gamma_2\circ\gamma_{(1,1)}^{-1}$ can be deformed to $\rho$.
If $\gamma_2\circ\gamma_{(1,1)}^{-1}$ can be deformed to $\rho$, we have
$C_2\cup_{\gamma_2^{-1}\circ\gamma_{(1,1)}} C_1\simeq
C_2\cup_{\rho} C_1$ and Proposition 3.4.1 will give the vanishing result.

The first approximation to our artificial transition map is constructed in
\cite{KoM}.
The domain is a neighborhood of the boundary of the cap $C_1$.
Consider the gluing data $(A_0,[\phi_x,F_x,y,(A_1,t_1),(A_2,t_2)])\in C_1$,
where $A_0$ is a connection on $Q_\alpha\times_{S^1}SO(3)$ which is
anti-self-dual
with respect to a metric $g_t$,
$\phi_x\in Q_\alpha|_x, F_x\in Fr(X)|_x,y\in D^4, x\in X$ and $(A_i,t_i)$ are
framed, charge one
instantons on $S^4$.
First, we compose the splicing map $\gamma_{(1,1)}'$ with the map $k$, this
gives an almost anti-self dual connection on $X$.
The frame $F_x$ identifies a ball around $x$, $D^4_x$ with $\bR^4$.
Parallel radial translation (with respect to the connection $A_0$)
of $\phi_x$ gives a trivialization of $Q_\alpha\times_{S^1}SO(3)$ over $B^4_x$.
If we write $P$ for the $SO(3)$ bundle obtained by gluing in the instantons,
this trivialization of $Q_\alpha\times_{S^1}SO(3)$ over $D^4_x$ gives one
of $P$ over an annulus $A_x$ whose outer boundary is equal to that of $B^4_x$.
This trivialization makes $\gamma_{(1,1)}'\circ k$ a connection on $\bR^4$,
framed at infinity.
Cutting off this connection, we can extend it to a connection $S^4$.
Keeping the radially parallel trivialization from the $Q_\alpha$ frame
gives a framed connection on $S^4$.  This gives a map we call $\rho'$.
We see $\rho'$ is $SO(4)$ equivariant with respect to the rotation action  and
the $SO(3)$ action.
The image of $\rho'$ is an almost anti-self dual connection on $S^4$, so
it can be perturbed to an anti-self dual one and then translated so it
is centered.  Both perturbation and translation commute with the rotation
action so this operation is $SO(4)$ and $SO(3)$ equivariant on $C_1^s(x)$.
Call this map $\rho(A_0,x)$.  Let $\rho_1$ be the extension of $\rho(A_0,x)$
to $\partial C_1$.

\

{\bf Remark}
There is a deformation given in \cite{KoM} between these transition maps,
but the author is uncomfortable with their argument.
They assert that the two gluing maps $\gamma_{(1,1)}$ and $\gamma_2\circ\rho$
are close enough to imply the existence of a deformation between them.
The best known
bounds on the size of the perturbation map are $L^2_1$ bounds in terms
of the dilation parameter of the $S^4$ connection being glued in.
Better bounds (i.e. $L^2_k$) would not seem possible as the cut-off function
in the splicing construction would not allow better bounds on $F^+_A$
than $L^2$.  The images of the maps $\gamma_{(1,1)}$ and $\gamma_2\circ\rho$
can be made close in the $L^2_1$ metric topology, but only by restricting
their domains to more concentrated connections.  In the $L^2$ metric topology,
the ends of the moduli space look like cones and we would be making the two
different embeddings close by pushing them down to the cone point.
It seems reasonable that the $L^2_1$ metric topology does not see this cone
structure (the $L^2_1$ metric is proportional the the conformally invariant
$L^4$
metric on the $S^4$ moduli space), but there is a great deal to be checked.
An explicit deformation between the two maps is presented here instead.

\

\proclaim{Proposition 3.4.2}
There is a deformation through embeddings between
the transition maps $\rho_1$ and $\gamma_2^{-1}\circ\gamma_{(1,1)}\circ k $.
\endproclaim
\demo{Proof}
Our goal is to show that the map $\gamma_{(1,1)}\circ k$ can be deformed,
through embeddings,
to $\gamma_2\circ\rho_1$.
We begin with the observation that the splicing map $\gamma_{(1,1)}'\circ k$
is equal to the splicing map $\gamma_2'\circ \rho'$.
There are two ways to perturb the images of these two spaces onto the
moduli space, arising from the two ways of creating right inverses
to the linearization of the anti-self dual equation.
Recall from \cite{DK}, Section 7.2 that a right inverse to $d^+_{A'}$
where $A'$ is in the image of $\gamma_{(1,1)}'\circ k$ can be constructed
from splicing together right inverses from the background connection
and the $S^4$ instantons.  The splicing construction depends on
the frames and cut-off functions being used and thus will be different
from the $\Sigma_{(1,1)}$ and $\Sigma_2$ construction.  Note that
although the $S^4$ connection in the image of $\rho'$ is not anti-self-dual,
it is sufficiently close to the anti-self-dual connections
that this construction can still be applied to it.
We thus have two different right inverses, $P_{(1,1)}$ and $P_2$ to $d^+_{A'}$
and thus two different perturbations of the image of
$\gamma_{(1,1)}'\circ k$ to the moduli space.
We construct a family of right inverses $P_t=(1-t)P_{(1,1)}+tP_2$ to $d^+_{A'}$
which will satisfy the necessary bounds (e.g. Lemmas 7.2.18 and 7.2.23 in
\cite{DK})
and thus define a family of perturbations of the splicing map to the moduli
space.
Our first deformation
of $\gamma_{(1,1)}\circ k$ is to change the perturbation
of $\gamma_{(1,1)}'\circ k$ to the moduli space by using this family of right
inverses.
Thus we have deformed the map $\gamma_{(1,1)}\circ k$ to the map
given by composing $\gamma_2'\circ\rho'$ with the perturbation to the moduli
space
defined by a right inverse constructed from the gluing data of the lower
stratum.
As noted above, although the image of $\rho'$ is not anti-self dual, it is
close
enough the the anti-self dual connections to have a right inverse with the
appropriate
bounds and the map $\gamma_2\circ\rho'$ can be defined by adding a perturbation
of
$\gamma_2'\circ\rho'$ to make the image anti-self dual.
Now, the map $\rho_1$ is a deformation of $\rho'$ (the perturbation is a
deformation).
We deform $\gamma_2'\circ\rho'$ to $\gamma_2'\circ\rho_1$ and
thus deform the maps to the moduli space $\gamma_2\circ\rho'$ to
$\gamma_s\circ\rho$.
We end up with the map $\gamma_2\circ\rho_1$.
\enddemo

Our new transition function $\rho_1$ does not yet satisfy the conditions of
Proposition
3.4.1 because it depends on the background connection $A_0$ and
the metric $g_t$.  But our family of background
connections and metrics is contractible, so the map $\rho_1$
can be deformed to $\rho_1([\alpha],x)$,
where $\alpha$ is the reducible connection.

\enddocument
\end